\documentclass[11pt]{amsart}

\usepackage[utf8]{inputenc} \usepackage[T1]{fontenc} \usepackage{lmodern}
\usepackage{pgf,csquotes}
\usepackage{epsfig,graphicx,color,amsmath,a4wide,amssymb,amsfonts,amsbsy,xy,latexsym,epsf,verbatim}
\usepackage{newtxtext}

\usepackage[
    backend=biber,
    style=alphabetic,
    maxcitenames=2, %
    maxbibnames=99, %
    giveninits=true, %
    sorting=nyt, %
    uniquename=false,
    uniquelist=false, %
    sortcites=true,
    hyperref=true]{biblatex}
\AtEveryBibitem{\clearfield{review} \clearfield{series} \clearfield{doi} \clearfield{isbn} \clearfield{issn} \clearfield{language} }

\DeclareFieldFormat{url}{Available at \url{#1}}
\makeatletter
\DeclareFieldFormat{eprint:arxiv}{%
  arXiv\addcolon\space
  \ifhyperref
    {\href{https://arxiv.org/\abx@arxivpath/#1}{%
       \nolinkurl{#1}%
       \iffieldundef{eprintclass}
         {}
         {\addspace\mkbibbrackets{\thefield{eprintclass}}}}}
    {\nolinkurl{#1}%
     \iffieldundef{eprintclass}
       {}
       {\addspace\mkbibbrackets{\thefield{eprintclass}}}}}
\makeatother

\definecolor{mylinkcolor}{rgb}{0.0,0.0,0.75}
\definecolor{myurlcolor}{rgb}{0.0,0.0,0.75}
\usepackage[
	colorlinks, urlcolor=myurlcolor,citecolor=myurlcolor,linkcolor=mylinkcolor,
        pdfusetitle,
]{hyperref}

\usepackage{algpseudocode}

\usepackage{tikz}
\usetikzlibrary{positioning}
\usetikzlibrary{cd}
\usetikzlibrary{arrows.meta}
\usetikzlibrary{circuits.ee.IEC}

\usepackage{cancel}

\usepackage[most]{tcolorbox}

\newtcolorbox{linebox}[2][]{
  enhanced,
  left=2pt, right=2pt, top=8pt, bottom=8pt,
  colback=white, colframe=black,
  boxrule=1pt, arc=3pt,
  #1,
  title=#2,
  coltitle=black, titlerule=0pt,
  fonttitle=\scriptsize\itshape,
  attach boxed title to bottom center={yshift=3mm, xshift=3cm},
  boxed title style={size=small, colback=white, colframe=white}
}

\usepackage{pgf,color}
\usepackage{booktabs}
\usepackage{multirow}

\definecolor{LightGrey}{rgb}{.85,.85,.85}
\definecolor{DarkGrey}{rgb}{.5,.5,.5}

\definecolor{Blue}{rgb}{.0,.0,0.9}
\definecolor{LightBlue1}{rgb}{.2,.4,0.9}
\definecolor{LightBlue2}{rgb}{.3,.5,0.9}
\definecolor{LightBlue3}{rgb}{.4,.6,0.9}
\definecolor{LightBlue4}{rgb}{.5,.7,.9}
\definecolor{LightBlue5}{rgb}{.6,.8,.9}
\definecolor{LightBlue6}{rgb}{.7,.9,.9}

\definecolor{Red}{rgb}{.9,.0,.0}
\definecolor{LightRed1}{rgb}{0.9,.2,.4}
\definecolor{LightRed2}{rgb}{0.9,.3,.5}
\definecolor{LightRed3}{rgb}{0.9,.4,.6}
\definecolor{LightRed4}{rgb}{.9,.5,.7}
\definecolor{LightRed5}{rgb}{.9,.6,.8}
\definecolor{LightRed6}{rgb}{.9,.7,.9}

\def\giota{{\langle \iota \rangle}}
\def\gt{{\langle t \rangle}}
\def\gtp{{\langle t' \rangle}}
\def\g2t{{\langle 2t \rangle}}
\def\gtp{{\langle t' \rangle}}
\def\gT{{\langle T \rangle}}

\xyoption{all}
\usepackage[english]{babel}

\def\ev{\mathop{\rm{ev}}\nolimits }

\def\KK{{\mathbf K}}
\def\PP{{\mathbf P}}
\def\RR{{\mathbf R}}
\def\Ks{{{\mathbf K}_s}}
\def\Kss{{{\mathbf K}^*_s}}
\def\LL{{\mathbf L  }}
\def\ZZ{{\mathbf Z}}
\def\agot{{\mathfrak a}}
\def\bgot{{\mathfrak b}}
\def\cgot{{\mathfrak c}}
\def\dgot{{\mathfrak d}}
\def\egot{{\mathfrak e}}
\def\fgot{{\mathfrak f}}
\def\ggot{{\mathfrak g}}
\def\hgot{{\mathfrak h}}
\def\igot{{\mathfrak i}}
\def\jgot{{\mathfrak j}}
\def\lgot{{\mathfrak l}}
\def\mgot{{\mathfrak m}}
\def\ngot{{\mathfrak n}}
\def\pgot{{\mathfrak p}}

\def\tgot{{\mathfrak t}}
\def\ugot{{\mathfrak u}}
\def\vgot{{\mathfrak v}}
\def\xgot{{\mathfrak x}}

\def\cC{{\mathcal C}}
\def\cD{{\mathcal D}}

\def\cF{{\mathcal F}}
\def\cG{{\mathcal G}}

\def\cL{{\mathcal L}}

\def\cQ{{\mathcal Q}}

\newtheorem{definition}{Definition}
\newtheorem{lemma}[definition]{Lemma}
\newtheorem{proposition}[definition]{Proposition}
\newtheorem{theorem}[definition]{Theorem}

\newtheorem{remark}[definition]{Remark}

\ExplSyntaxOn
\cs_new_eq:NN \calc \fp_eval:n
\ExplSyntaxOff

\newenvironment{myproof}[1][\myproofname]{\par
  \normalfont \topsep6pt\relax
  \trivlist
\item[\hskip\labelsep
  \itshape
  #1.]\ignorespaces
}{%
  \endtrivlist\hfill$\square$
}
\providecommand{\myproofname}{Proof}

\author[Couveignes]{Jean-Marc Couveignes}
\address{Jean-Marc Couveignes, Univ. Bordeaux, CNRS, INRIA,
  Bordeaux-INP, IMB, UMR 5251, F-33400 Talence, France.}
\email{jean-marc.couveignes@u-bordeaux.fr}

\author[Lercier]{Reynald Lercier}
\address{%
  Reynald Lercier,
  DGA \& Univ. Rennes, %
  CNRS, IRMAR - UMR 6625, F-35000
  Rennes, %
  France. %
}
\email{reynald.lercier@univ-rennes.fr}

\title{Elliptic butterflies}

\usepackage{caption}
\usepackage{svg}

\begin{document}

\vspace*{-1.8cm}

\begin{abstract}
  We study natural evaluation and interpolation problems
  for elliptic functions and prove that they allow a recursive
  treatment using a variant of classical  butterflies  first
  introduced by Gauss. We deduce the existence of
  straight-line programs with
  complexity scaling with $d\log(d)$ for these problems and present applications
  to finite field  arithmetic, coding theory and cryptography.
\end{abstract}
\dedicatory{In memory of Tony Ezome}
\maketitle

\medskip

\begin{figure} [h!]
  \centering
  \def\bfratio{0.75}
  \begin{tcolorbox}[
    enhanced, notitle,
    colback=white, colframe=white,
    boxrule=0pt, arc=9pt,
    watermark graphics=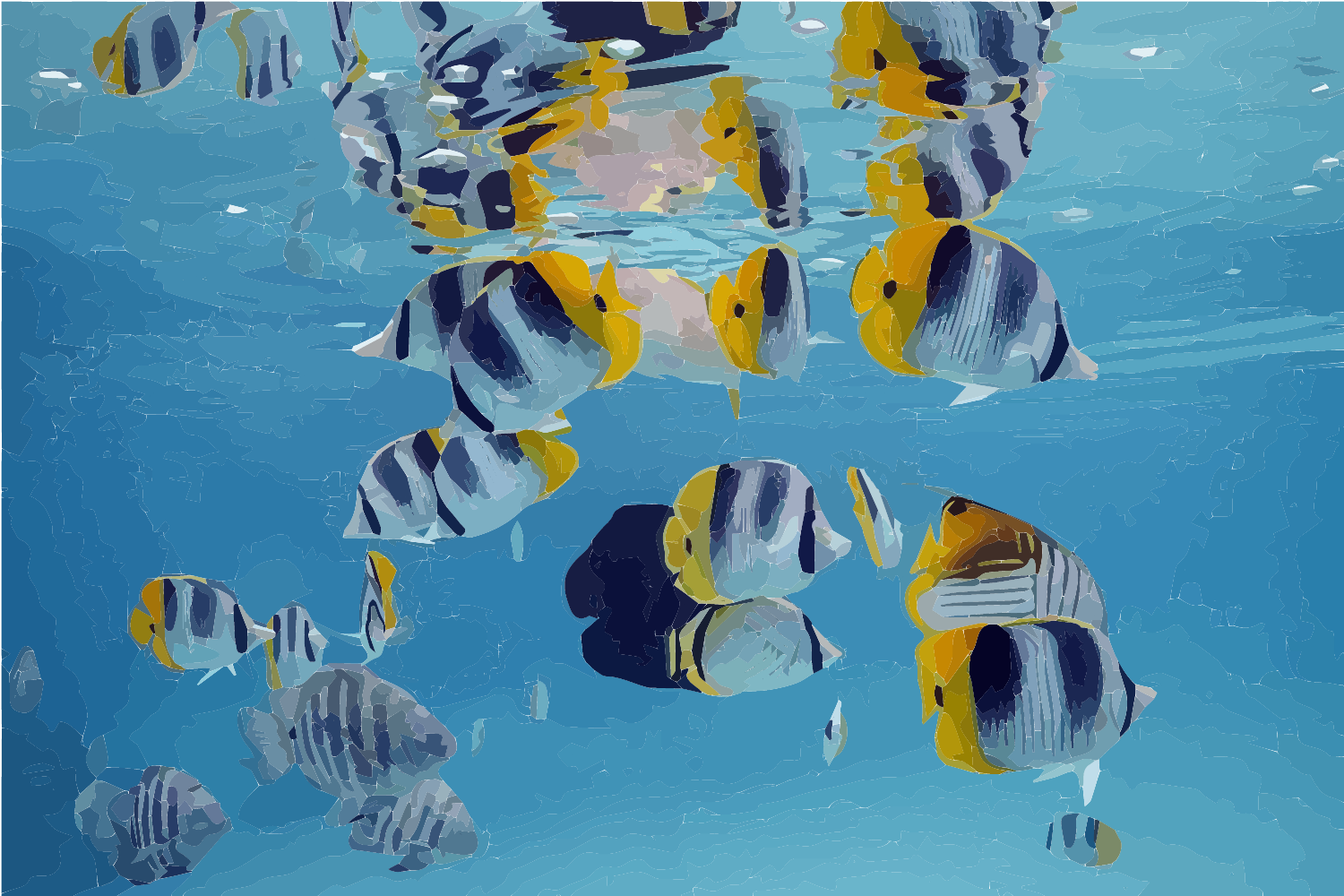,
    watermark overzoom=1.0,
    height=\calc{\bfratio*6}cm, width=\calc{\bfratio*9}cm]
  \end{tcolorbox}
  \caption*{\hfill \small Butterflyfish \tiny (freely adapted
    from~\cite{wiki})}
\end{figure}

  \tableofcontents
\newpage

\section{Introduction}\label{sec:intro}

Let $\KK$ be a field and let  $X$ be  a smooth, projective,  absolutely integral curve over $\KK$.
Let $D$ be a divisor on $X$. We denote by $\cL(D)$ the associated linear space and call $d$
its dimension.
Let $b_0$, $b_1$, \ldots, $b_{d -1}$ be $\KK$-points on $X$ and  not in the support of $D$.
We consider the evaluation map
\begin{equation}\label{eq:ev}
  \begin{tikzcd}[
    row sep = 0pt,
    /tikz/column 1/.style={column sep=0em},
    /tikz/column 2/.append style={anchor=base east},
    /tikz/column 3/.append style={anchor=base west} ]
    \ev : & \cL (D) \arrow[r]  & \KK^d\,, \\
    & f \arrow[r, mapsto] & f(b_l)_{0\leqslant l \leqslant d -1}\,.
  \end{tikzcd}
\end{equation}
which we assume to be invertible.
We are interested in the complexity of { evaluation} (evaluating the map $\ev$)
and
{ interpolation} (evaluating the reciprocal map
$\ev^{-1}$).
 Elements in $\KK^d$ are represented by their entries.
Elements in $\cL (D)$ are represented by their coordinates
in a basis $u=(u_0, u_1, \ldots, u_{d-1})$ that must  be specified. It is traditional
to use bases satisfying one or several among the following properties:
\begin{enumerate}
\item the sum of the degrees of the $u_l$ is minimal among all bases,
\item the $u_l$ are eigenvalues of some automorphism of $X$,
\item the $u_l$ are permuted by some automorphism of $X$.
\end{enumerate}
One  looks for straight-line  programs that solve the evaluation and interpolation problems using
additions and multiplications by  precomputed
constants in $\KK$.
Such programs contain no loops or
conditionals, and consist of a sequence of operations
involving the inputs, the constants, and
previously computed elements. The complexity
is then measured by the
number of operations.

The classical and most important case is when $X=\PP^1$ and $D=(d-1) [\infty]$
and $\cL(D)$ is the space of polynomials of degree
$\leqslant d -1$ in one variable $x$. In this situation, and assuming
that $\KK$ contains
a primitive $d$-th root
of unity, the monomial basis $(1,x , x^2, \ldots, x^{d-1})$ satisfies
the first two conditions above.
If further $d$ is a power of two, which implies that
the characteristic of $\KK$ is not two,  then
the complexities of evaluation and interpolation are bounded by \begin{equation}\label{eq:3/2}(3/2)d\log (d) +1.\end{equation} This is achieved using the well-known
fast Fourier transform or FFT, first discovered by Gauss.
See~\cite[Theorem 2.6]{BCS} and the historical account in~\cite{HJB}.

Bostan and Schost give in~\cite{BS}  a range of nice upper bounds for the complexity of evaluation
and interpolation  when the basis is the monomial basis or the Newton basis, and the evaluation set is
an arithmetic sequence, a geometric sequence, or any set. They
make no  assumption on $\KK$.  For  example, without making any restriction on $\KK$ nor on the evaluation
set, the complexity
of evaluation and interpolation in the monomial basis is  bounded by an absolute  constant times
$M(d) \log (d)$ where $M(d)$ is the complexity of multiplying two polynomials of degree $\leqslant d$
with coefficients in $\KK$. Since the latter is bounded by a constant times
$d\log (d)  \log (\log (d))$ operations in $\KK$,
according to a result of Sch\"onhage and Strassen~\cite[Theorem 2.13]{BCS}, we end
up with a general bound of a constant times \begin{equation}\label{eq:general}d
  \log^2 (d) \log (\log (d))\end{equation} for evaluation and interpolation
in the monomial basis, at a general set of points.

Comparing Equations (\ref{eq:general}) and (\ref{eq:3/2}) we see that we save
a factor $\log (d) \log (\log (d))$ when $d$ is a power of  two,
$\KK$ contains a primitive $d$-th root of unity,  $\omega_d$,   and
the evaluation set is the set of all  $d$-th roots of $1$.
We notice that in the specific case
of evaluation, resp. interpolation, in the monomial
basis over
a geometric progression of size~$d$, Bostan and Schost~\cite[Prop.~5]{BS} also
give an algorithm whose complexity decreases to
$M(d) + O(d)$,  resp. $2M(d)+O(d)$. This still
differs from Equation (\ref{eq:3/2}) by a factor of $\log (\log (d))$.
This gain is not marginal as is demonstrated by the many implementations and applications of FFT.
Indeed general evaluation and interpolation methods often use  specific ones as subroutines.

The key idea behind FFT is to use the involution $x\mapsto -x$  to decompose the  polynomial to be evaluated
$P(x)=\sum_{0\leqslant l\leqslant d-1} p_l\, x^l$ as a sum \[P(x)=P^+(x)+P^-(x)\]
where \[P^+(x) = \sum_{0\leqslant 2l\leqslant d-1} p_{2l}\, x^{2l} \text{\,\,\,
and \,\,\,} P^-(x)=\sum_{0\leqslant 2l+1\leqslant d-1} p_{2l+1} \,x^{2l+1}  = x\sum_{0\leqslant 2l+1\leqslant d-1} p_{2l+1} \, x^{2l}\] are the even and odd parts of $P(x)$. This decomposition
reduces the evaluation of $P(x)$ to two similar problems of halved size, allowing a recursive approach that
is classically illustrated using diagrams called butterflies (see Figure~\ref{fig:butterfly}).

\begin{figure}[htbp]
  \resizebox{!}{5cm}
  {%
    \begin{tikzpicture}[
      xscale=1.5, thick, node distance=.5cm, circuit ee IEC,
      font=\footnotesize,
      box/.style={
        draw, align=center, shape=rectangle, minimum width=3cm,
        minimum height=3cm, append after command={
          \foreach \side in {east,west} {
            \foreach \i in {1,...,#1} {
              (\tikzlastnode.north \side)
              edge[draw=none, line to] coordinate[pos=(\i-.5)/(#1)]
              (\tikzlastnode-\i-\side) (\tikzlastnode.south \side) }
          }
        }
      }
      ]

      \node[box=4] (box-t) {\footnotesize 4-point FFT};
      \node[box=4, below=of box-t] (box-b) {\footnotesize 4-point FFT};

      \foreach \s[count=\i] in {0,2,4,6}
      \path (box-t-\i-west) edge node[at end, left]{$p_{\s}$} ++(left:.5);
      \foreach \s[count=\i] in {1,3,5,7}
      \path (box-b-\i-west) edge node[at end, left]{$p_{\s}$} ++(left:.5);

      \foreach \b/\s[count=\k] in {t/+, b/-}
      \foreach \i[evaluate={\j=int(\i-1)}, evaluate={\J=int(ifthenelse(\k==2,\j+4,\j))}] in {1,...,4}
      \node [contact] (conn-\b-\i) at ([shift=(right:1.5)]
      box-\b-\i-east) {} edge node[above] {${P^{\s}(\omega_4^{\j})}$} (box-\b-\i-east)
      node [contact, label=right:{${P(\omega_8^{\J})}$}] (conn-\b-\i') at
      ([shift=(right:5)] box-\b-\i-east) {};

      \begin{scope}[every info/.append style={font=\tiny, inner sep=+.1pt}]

        \foreach \i[evaluate={\j=int(\i-1)},evaluate={\J=int(\i+3)}] in {1,...,4}
        \path (conn-t-\i)
        edge[current direction={pos=.27}] (conn-t-\i')
        edge[current direction={pos=.1}] (conn-b-\i') (conn-b-\i)
        edge[current direction={pos=.85, info={$\omega_8^\J$}}] (conn-b-\i')
        edge[current direction={pos=.88, info={$\omega_8^\j$}}] (conn-t-\i');
      \end{scope}
    \end{tikzpicture}
  }
  \caption{An 8-point FFT butterfly}
  \label{fig:butterfly}
\end{figure}

In this work we define  and study a family of evaluation and interpolation problems
on an elliptic curve $E$ having a $\KK$-rational point $t$ of order $d$,
a power of two. We show that there are natural
analogues of butterflies in this context. We replace   roots of
unity  by multiples of the point $t$
and the associated translations,  while polynomials are replaced by functions
in the space $\cL(\gt)$ equipped with a basis $u$ of functions of degree $\leqslant 2$.
In particular the role of the  involution $x\mapsto -x$ is now played
by the  translation by $T=(d/2)t$, a point of order   two.
This results in complexity bounds
$\cQ  d  \log (d)$
for the evaluation and   interpolation problems\footnote{Note that everywhere in this article, the notation $\cQ$
stands for a positive absolute constant. Any sentence containing
this symbol becomes true if the symbol is replaced in every occurrence by some large enough real number.}.

To be complete we should say that evaluation and interpolation
are not quite enough for our purposes. We  need another
linear map on functions that we call { reduction}
and which
we  now define in a general context. We now have three divisors $D$, $E$ and $B$ on the curve
$X$.
We denote by $d$ the dimension of $\cL(D)$.
We assume that $B$ is effective
of  degree  $d$. We  assume
that $D$ and $B$ have  disjoint  support.
We denote by $\LL$ the residue ring at $B$. This is a $\KK$-algebra
of dimension $d$. We assume that the residue map
\begin{equation*}
  \begin{tikzcd}[
    row sep = 0pt,
    /tikz/column 1/.style={column sep=0em},
    /tikz/column 2/.append style={anchor=base east},
    /tikz/column 3/.append style={anchor=base west} ]
    \ev_D^B : & \cL (D) \arrow[r]  & \LL\,,\\
    & f \arrow[r, mapsto] & f \bmod B\,.
  \end{tikzcd}
\end{equation*}
is invertible.
We  also assume
that $E$ and $B$ have  disjoint  support and define
\begin{equation*}
  \begin{tikzcd}[
    row sep = 0pt,
    /tikz/column 1/.style={column sep=0em},
    /tikz/column 2/.append style={anchor=base east},
    /tikz/column 3/.append style={anchor=base west} ]
    \ev_E^B : & \cL (E) \arrow[r]  & \LL\,,\\
    & f \arrow[r, mapsto] & f \bmod B\,.
  \end{tikzcd}
\end{equation*}
the residue map.
For every function $F$ in $\cL(E)$ there exists a unique
function $f$ in $\cL(D)$ such that \[F\equiv f\pmod B
  \text{\,\,\,
that is \,\,\,} \ev_E^B(F)=\ev_D^B(f)\in \LL.\]
We choose a basis for $\cL(D)$ and a basis for $\cL(E)$.
The reduction problem takes as input the coordinates
of $F$ and returns the coordinates of $f$.
We will explain how to efficiently solve
a family of reduction problems on elliptic curves.
In the situation we are interested in, the divisor $E$ is $2D$ because we need
to reduce the product of two functions in $\cL(D)$.

We will  present three applications of these fast evaluation,
 interpolation, and reduction  algorithms. The first one is the construction
of normal bases for  extensions of finite fields of degree a power of two,
allowing  particularly fast multiplication.
The  second one is the description of  $[d,d/2,d/2+1]$-error correcting codes
 that can be encoded and checked in time $\cQ d \log (d)$.
The third one is a discussion on the relevance of using residue rings on elliptic
   curves in the context of LWE cryptography.
\medskip

 To conclude this introduction
 we recall  that there is a long tradition of interpolating on algebraic curves especially over finite fields.
 Among the major results, one may mention
linear upper bounds  on the bilinear complexity of multiplication in finite field extensions
by  Chudnovsky and Chudnovsky ~\cite{chu},
Shparlinski, Tsfasman, Vl\u{a}du\c{t}~\cite{STV},
Shokrollahi~\cite{SH}, Ballet  and Rolland~\cite{BR,BAL},
Chaumine~\cite{CHA},  Randriambololona~\cite{RAN}
and others. Another major achievement was the construction of excellent codes  by Goppa, Tsfasman, Vl\u{a}du\c{t}, Zink, Ihara, Garc\'{\i}a, Stichtenoth ~\cite{gop1, gop2, tvz,garsti,ihara}.
Chudnovsky  and Chudnovsky express in~\cite[Section 6]{chchch} the intuition that elliptic curves with smooth order could be of some use to accelerate some specific interpolation  problems in genus $0$.
Using similar ideas,
Ben-Sasson,  Carmon,  Kopparty,  and Levit show in~\cite{bckl}
how to quickly compute the values of a polynomial at a set $S'$ from its value at another set
$S$ when $S$ and $S'$ are specific subsets of finite fields
{\it i.e.}  $x$-coordinates of well chosen points on a well chosen elliptic curve.
Starting from Chudnovsky's intuition we prefer to stay in the elliptic world  where all these ideas and objects originate from, and develop an explicit analogue
of butterflies in this context. We also prefer to invoke evaluation and interpolation rather than   Fourier transform,  because in our humble opinion, a Fourier transform
is something different. See \text{e.g.}~\cite[Section 13.5]{BCS}
or~\cite[Section 6.1]{cg}.
\smallskip

Finally, we would like to stress that we do not address the problem
of polynomial multiplication in this work. We stick
to elliptic functions.
In principle, one could disguise two polynomials as elliptic functions, then
multiply them out
using elliptic butterflies. The cost of translation,
however, would be $O(d\log^2(d))$ just as in~\cite{bckl}.
This would not be competitive with Sch\"onhage-Strassen method~\cite[Theorem 2.13]{BCS}.
We consider its $\cQ d\log (d) \log (\log (d))$ complexity
as the state of the art in polynomial multiplication
and use this upper bound whenever needed in this article,
despite the existence of an algorithm
by Harvey and van der Hoeven~\cite{HvdH}
that multiplies polynomials with coefficients
in a finite field with
heuristic complexity $\cQ d\log (d)$. Consequently,
all the complexity estimates given in this article are fully
unconditional.
\medskip

\noindent \textbf{Acknowledgements.}
We thank the anonymous referees for a careful reading of the manuscript
and for helpful suggestions. We are also grateful to Grégoire Lecerf
for valuable comments.

 \section{Summary and organization of the article}

 In this article, we consider  specific evaluation and interpolation
 problems in genus one.
 We consider an elliptic curve $E$ over a field $\KK$ of odd characteristic.
 We let $d$ be a power of two and
 we assume that  there is a point $t$ of order $d$ in $E(\KK)$.
 We let \[D=\sum_{0\leqslant l\leqslant d-1}[l t]\] be the degree $d$ divisor
 associated to the group $\gt$ generated by $t$. We take $b$ to be a $\KK$-rational point on $E$
 such that $db$ is not equal to the point at infinity. We set \[b_l=b+l t \text{\,\,\, for \,\,\,} l\in [0, d-1]\] and we check that none of these $d$
 points is in $\gt$.
 The evaluation map $\ev$ in Equation (\ref{eq:ev}) is well defined and invertible.

 To complete the description of the interpolation problem we are interested in, we need interesting
 bases for $\cL (D)$. These bases are made of functions of degree $\leqslant 2$. In Section \ref{sec:ef} we  recall
 nice properties of degree two functions on an elliptic curve. We use these functions in Section
 \ref{sec:ebfu} to define two nice bases for $\cL(D)$. One of these bases  satisfies the first condition stated
 in the introduction. The sum of the degrees of its elements is $2d-2$. The other basis  satisfies the third condition.
 It is invariant by translation by $t$. And the sum of the degrees of its elements is $2d$.
 We show that base change between these two bases requires no more than
 $\cQ  d$ operations in $\KK$. This means that these two
 bases are essentially the same basis from the point
 of view of complexity theory.

 The automorphism $x\mapsto -x$ is crucial to define classical
butterflies. Its  elliptic counterpart is  translation
by $T$ where $T=(d/2)t$ is the unique point of order two in $\gt$.
We study this translation and the corresponding quotient $E'=E/\gT$  in Section \ref{sec:2}. This enables us to describe
the recursive evaluation algorithm in Section \ref{sec:but}. An important
point is that we use a degree two function
$\theta$ that is odd for translation by $T$. We decompose $f\in \cL(D)$ as \[f=f^++f^-\] where $f^+$ is even
and $f^-$ is odd. We then set \[f^0=f^+ \text{\,\,\,  and \,\,\,}
f^1=\theta  f^- \] and we notice that both are functions on the quotient $E'$.
In the classical setting the role of $\theta$ is played by the variable $x$  itself.

The interpolation algorithm is already a bit more tricky than evaluation. It is described in Section \ref{sec:int}.  A difficulty is to solve a linear system of a very special and classical form: a cyclic bidiagonal system. We recall in Section \ref{sec:bidia}
what is needed to solve efficiently such a system in our context.

If we multiply two functions $f_1$ and $f_2$
in $\cL (D)$ we obtain a function $F$ in $\cL(2D)$. Given such  a function $F$, there is a unique function $f$ in $\cL(D)$
such that  $f(b_l)=F(b_l)$ for every $l\in [0,d-1]$. Finding $f$ once given $F$ is a reduction problem.
In Section \ref{sec:redu} we define and solve
a natural generalization of this problem. The resulting
efficient multiplication algorithm in the residue ring
at $B$ is presented in Section \ref{sec:resid}.
Section \ref{sec:impl-exper-results} documents
a public implementation of the algorithms presented in this article and
provides experimental results.

A first application, presented  in Section \ref{sec:ebff},
is the construction of
normal bases for finite field extensions of degree a power of two,
allowing multiplication in time $\cQ d \log (d)$. We build on  the construction in~\cite{CLEP} of an equivariant
version of Chudnovsky's  interpolation method~\cite{chu},  and we use our complexity estimates for evaluation, interpolation and reduction
in this context.

In Section \ref{sec:gop}, we construct Goppa\footnote{In this
  article, we use the term Goppa code for any code
  obtained by evaluation of algebraic functions on a given algebraic
  curve at a collection of points
on that curve.} codes in genus one obtained by restricting the evaluation map
$\ev$ to the subspace of functions that are invariant by the elliptic involution $P\mapsto -P$
on $E$. The resulting  code is  a genus zero Goppa code. The advantage of seeing it as a genus one code is that we benefit from
the acceleration due to the presence of the automorphism group $\gt$. This produces MDS $[d,d/2,d/2+1]$-codes that can be encoded and checked
in time $\cQ d \log (d)$ and decoded up to half the minimum distance at the expense of
$\cQ d \log^2 (d)\log (\log (d))$ operations and comparisons
in $\KK$. The only constraint for the existence of such codes is that the size of $\KK$
be large enough with respect to $d$. This extends  Lacan and Soro's result~\cite{laso} for Fermat fields.

The purpose of the final Section \ref{sec:ellipt-lwe-crypt} is to discuss
the advantages and  limitations of using elliptic residue rings
in the context of LWE-cryptography.

All the straight-line
programs presented in this work are summarized
in pseudo-code right after their  algorithmic description.
An  implementation of these algorithms  and the underlying
data structures
are introduced in Section \ref{sec:data} and Section~\ref{sec:impl-exper-results}.

\section{Algorithms and data structures}\label{sec:data}

The algorithms presented in this article are designed to manipulate elements
in vector spaces. Vectors are typically  coordinate vectors
of functions in
well-chosen
bases of Riemann–Roch spaces, or  evaluation vectors of functions
at well chosen points. We therefore use
the classical $d$-dimensional vector notation $\vec{a}$ and
$(a_l)_{l\in \ZZ/d\ZZ}$ for inputs, outputs, variables and constants.

For ease of reading, our pseudo-code descriptions of these algorithms are
given in a recursive form. However, like the classical butterfly formulation
of FFT algorithms \emph{à la Cooley–Tukey}, they can be very easily unrolled
to yield iterative functions. In fact, our algorithms are all straight-line
programs.
\medskip

Beyond the usual operations of addition and scalar multiplication, denoted by
``$+$'' and “$\cdot$”, we introduce the following notations for other operations
on vectors.
\begin{itemize}
\item We denote by $\langle \vec{a},\, \vec{b}\rangle$ the scalar product of
  $\vec{a}$ and $\vec{b}$.
\item We denote by $\sigma(\vec{a})$ the   cyclic rotation,
\textit{i.e.} the
  vector $(a_{\sigma(l)})_{l\in \ZZ/d\ZZ}$, where $\sigma(l) = l-1$.
\item We denote by $e + \vec{a}$ the  addition by a constant, \textit{i.e.}
  the
  vector $(e+a_l)_{l\in \ZZ/d\ZZ}$.
\item We denote by $\vec{a}\,\oplus\, \vec{b} $ the direct sum, \textit{i.e.}
  the concatenation of $\vec{a}$ and $\vec{b}$.
\item We denote by $\vec{a}\,\diamond\, \vec{b}$ the component-wise product of
  $\vec{a}$ and $\vec{b}$.
\item We denote by $ \vec{a}\,^{\mathrm{inv}} = (a_l^{-1})_{l\in
    \ZZ/d\ZZ}$ the inverse vector of $\vec{a}$ for the product $\diamond$\,.
\end{itemize}

\section{Elliptic functions}\label{sec:ef}

We recall a few formulae from~\cite{CLEP} regarding low-degree functions
on elliptic curves.
We let  $\KK$ be a field
 and  $E$  the elliptic curve over $\KK$
defined by the Weierstrass equation
\begin{displaymath}
Y^2Z+a_1XYZ+a_3YZ^2=X^3+a_2X^2Z+a_4XZ^2+a_6Z^3\,.
\end{displaymath}
Call $O=(0,1,0)$ the point at infinity.
We set $x=X/Z$,
$y=Y/Z$ and $z=-x/y=-X/Y$. The Taylor expansions
of $x$ and $y$ at $O$ in the local parameter $z$ are
\begin{eqnarray*}
x&=&\frac{1}{z^2}-\frac{a_1}{z}-a_2-a_3z+O(z^2)\,,\\
y&=&-\frac{1}{z^3}+\frac{a_1}{z^2}+\frac{a_2}{z}+a_3+O(z)\,.
\end{eqnarray*}
If $A$ is a point on $E$, we denote by
$\tau_A$ the translation by $A$.
We denote  by \[z_A = z \circ \tau_{-A}\] the composition
of $z$ with the translation by $-A$. This is a local parameter
at $A$.
We define  $x_A$ and $y_A$ in a similar way.
If  $A$ and $B$ are two distinct  points on $E$, we denote
by $u_{A,B}$ the function  on $E$ defined as
\begin{equation*}\label{eq:slope}
u_{A,B}=\frac{y_A-y(A-B)}{x_A-x(A-B)}\,.
\end{equation*}
It has  polar divisor $[A]+[B]$.
 It is invariant by the
involution exchanging $A$ and $B$,
$$u_{A,B}(A+B-P)=u_{A,B}(P)\,.$$
If $A$, $B$ and $C$ are three  points on $E$
we denote by
$\Gamma(A,B,C)$ the slope of the secant (resp. tangent) to $E$  going through  $C-A$ and $A-B$.
As a function on $E^3$ it is well defined for any three points $A$, $B$, $C$ such
that $\# \{A,B,C\}\geqslant 2$.
 The  Taylor expansions of $u_{A,B}$ at $A$ and $B$ are
\begin{eqnarray}
\label{eq:TuAB} u_{A,B}&=&-\frac{1}{z_A}-x_{A}(B)z_A+(y_{A}(B)+a_3)z_A^2+O(z_A^3)\\
&=&
\nonumber \frac{1}{z_B}-a_1+x_{A}(B)z_B+(y_{A}(B)+a_1x_A(B))z_B^2+O(z_B^3).
\end{eqnarray}
We deduce \begin{eqnarray}\label{eq:TuAB'}
  \nonumber  u_{B,A}&=&-u_{A,B}-a_1,\\
             [0.2em]
  u_{A,B}+u_{B,C}&=&u_{A,C}+\Gamma(A,B,C),\label{eq:somme2}\\[0.2em]
\nonumber \Gamma(A,B,C)&=&u_{B,C}(A)=u_{C,A}(B)=u_{A,B}(C)=-u_{B,A}(C)-a_1,\\[0.2em]
\nonumber u_{B,C}&=&u_{B,C}(A)-(x_{A}(C)-x_{A}(B))z_A+(y_{A}(C)-y_{A}(B))z_A^2+ O(z_A^3),\\[0.2em]
\nonumber   u_{A,B}u_{A,C}&=&x_A+\Gamma(A,B,C)u_{A,C}+\Gamma(A,C,B)u_{A,B} +a_2
+x_{A}(B)+x_{A}(C),\\[0.2em]
\nonumber u_{A,B}^2&=&x_A+x_B-a_1u_{A,B}+x_A(B)+a_2.
\end{eqnarray}

\section{Bases for elliptic  function field extensions}\label{sec:ebfu}

Let $\KK$ be a field with odd characteristic and let $E$
be an elliptic curve over $\KK$.
Let $t$ be a point of order $d=2^\delta$ in $E(\KK)$.
We assume that $\delta\geqslant 1$.
The group $\gt$ generated by $t$
can be seen as an effective divisor of degree
$d$ on $E$. We let $\cL (\gt)$ be the associated $\KK$-linear space.
We
construct two bases for $\cL (\gt)$ consisting of functions
of degree $\leqslant 2$.
We start with the
following lemma from~\cite[Lemma 4]{CLEP}.

\begin{lemma}
The sum $\sum_{l\in \ZZ/d\ZZ} u_{lt,(l+1)t} $
is a constant $\agot \in \KK$.
\end{lemma}

\begin{myproof}
The sum $\sum_{l\in \ZZ/d\ZZ} u_{lt,(l+1)t} $ is invariant
by translations in $\gt$. So it can be seen
as a function on $E/\gt$. As such, it has no more than one
pole. Therefore it is constant. Repeated use
of Equation (\ref{eq:somme2}) shows that
\[\agot= -a_1+\Gamma(O,t,2t)+\Gamma(O,2t,3t)+\Gamma(O,3t,4t)+\dots
+\Gamma(O,(d-2)t,(d-1)t).\]
\end{myproof}

For $l$ in $\ZZ/d\ZZ$ we set \[u_l=u_{lt,(l+1)t}+(1-\agot)/d\]
and check that
\begin{equation}\label{eq:change}\sum_{l\in \ZZ/d\ZZ}u_l=1 \in \KK.\end{equation}
We let  $(v_l)_{l \in \ZZ/d\ZZ}$ be  the system defined by
\[v_0=1\text{ and }v_l=u_{O,lt}\text{ for } l\not = 0\bmod d.\]
Examination of poles shows that both
$(u_l)_{l \in \ZZ/d\ZZ}$ and $(v_l)_{l \in \ZZ/d\ZZ}$
are bases of $\cL (\gt)$. For every integer $l$ such that
$0\leqslant l\leqslant d-1$ there exists a constant $\agot_l$
such that
\begin{equation}\label{eq:change2}
  v_{l} = \sum_{m=0}^{l-1}u_m +\agot_l.\end{equation}
This is proved using
Equation (\ref{eq:somme2}) repeatedly.
We compute $\agot_0=1$,
$\agot_1=(\agot-1)/d$ and for $2\leqslant l\leqslant d-1$
\begin{equation}\label{eq:agot}
  \agot_l= l(\agot - 1 )/d- \Gamma(O,t,2t)-\Gamma(O,2t,3t)-\Gamma(O,3t,4t)-\dots
  -\Gamma(O,(l-1)t,lt).
\end{equation}

 In what follows, we denote by
\textsc{UtoV\_BaseChange} and \textsc{VtoU\_BaseChange}  the subroutines that
implement, respectively, the base change from
$(u_l)_{l \in \ZZ/d\ZZ}$ and $(v_l)_{l \in \ZZ/d\ZZ}$ and \textit{vice versa}.
From Equations (\ref{eq:change}) and
(\ref{eq:change2}), we deduce that these routines require no more than
$\cQ d$ additions and multiplications by a constant  in $\KK$.

\section{Cyclic bidiagonal linear systems}\label{sec:bidia}

We now consider linear sparse
systems with a very special
structure, called cyclic
bidiagonal.  These systems arise naturally in Section~\ref{sec:but}, and the
ability to invert them efficiently is an essential ingredient for
Section~\ref{sec:int}.\smallskip

Let $d\geqslant 2$  be an integer. Let $\RR$  be a commutative local
ring.
Let $b_0$, $b_1$, \ldots, $b_{d-1}$
and $c_0$, $c_1$, \ldots, $c_{d-1}$
be scalars in $\RR$. Let $M$ be the matrix

\smallskip

\[M=\begin{pmatrix}
    b_{0} & 0 & \cdots & 0 & c_{0}\\
    c_1 & b_1 & 0 & \cdots & 0\\
    0  & c_2 & b_2 & \ddots & \vdots\\
    \vdots & \ddots & \ddots & \ddots & 0 \\
    0 &\cdots & 0 & c_{d-1} & b_{d-1}
  \end{pmatrix}.\]

\smallskip

The determinant
of $M$ is \[\det (M) = \prod_{i=0}^{d-1} b_i-
  (-1)^{d}\prod_{i=0}^{d-1} c_i.\]
We  assume that it is   invertible in $\RR$.  The ring $\RR$
being local, this implies that either every $b_i$
is invertible  or every $c_i$ is invertible.
Let $s_0$, $s_1$, \ldots, $s_{d-1}$
be scalars in $\RR$.
We solve the system

\smallskip

\[
  \begin{pmatrix}
    b_{0} & 0 & \cdots & 0 & c_{0}\\
    c_1 & b_1 & 0 & \cdots & 0\\
    0  & c_2 & b_2 & \ddots & \vdots\\
    \vdots & \ddots & \ddots & \ddots & 0 \\
    0 &\cdots & 0 & c_{d-1} & b_{d-1}
  \end{pmatrix}
  \begin{pmatrix}t_0\\t_1\\\vdots \\t_{d-2}\\t_{d-1}\end{pmatrix}
  =  \begin{pmatrix}s_0\\s_1\\\vdots \\s_{d-2}\\s_{d-1}\end{pmatrix}\]

\medskip

\noindent
in the $d$ unknowns  $t_0$, $t_1$, \ldots, $t_{d-1}$.
We eliminate $t_0$ from the first equation  using the second
equation. We then
eliminate
$t_1$ from the resulting equation using the
third equation. And so on. We find
\[\det (M) \cdot  t_{d-1}= -
  \sum_{0\leqslant k \leqslant d-1}(-1)^{k}
  \left( \prod_{0\leqslant l\leqslant k-1} b_l \right) \cdot
  s_k \cdot \left(\prod_{k+1\leqslant l\leqslant d-1}c_l\right).\]

\smallskip
\noindent
Assume that all the $b_i$ are invertible.
We compute $t_0=(s_0-c_0t_{d-1})/b_0$,
$t_1=(s_1-c_1t_0)/b_1$, \dots,
$t_{d-2}=(s_{d-2}-c_{d-2}t_{d-3})/b_{d-2}$.
Otherwise we know that all the $c_i$ are invertible. We then compute
$t_{d-2}=(s_{d-1}-b_{d-1}t_{d-1})/c_{d-1}$,
$t_{d-3}=(s_{d-2}-b_{d-2}t_{d-2})/c_{d-2}$,
\dots,
$t_1=(s_2-b_2t_2)/c_2$,
$t_0=(s_1-b_1t_1)/c_1$ .\smallskip

Thus, both evaluating and inverting a cyclic bidiagonal system of dimension
$d$ require only $\mathcal{O}(d)$ operations. In the sequel, we denote by
\textsc{BiDiagonal\_Evaluate} and \textsc{BiDiagonal\_Invert} the routines that
implement, respectively, the direct matrix–vector product and this elimination
procedure.

\section{A Vélu  isogeny of degree two}\label{sec:2}

Under the  hypotheses and  notation at the beginning of Section
\ref{sec:ebfu} we call \[T=2^{\delta-1}t\] the point of order two in the group
$\gt$. Following Vélu ~\cite{JV,JV2}
we write  $T=(x(T),y(T))$ and set
\begin{align*}
        \begin{split}
          b_2&=a_1^2+4a_2,\\
          b_4&=a_1a_3+2a_4,\\
          b_6&=a_3^2+4a_6,\\
          w_4&=3x(T)^2+2a_2x(T)+a_4-a_1y(T),\\
          w_6&= 7x(T)^3 + (2a_2 + b_2)x(T)^2 + (a_4 + (-y(T)a_1 + 2b_4))x(T) + b_6,\\
          A_4&=a_4-5w_4,\\
          A_6&=a_6-b_2w_4-7w_6
\end{split}\end{align*}
and
\begin{align}
        \begin{split}\label{eq:vel1}
          x'&=x+x_T-x(T),\\
          y'&=y+y_T-y(T),\\
          z'&=-x'/y'.
          \end{split}\end{align}
We denote by
\begin{equation*}
  \begin{tikzcd}[
    row sep = 0pt,
    /tikz/column 1/.style={column sep=0em},
    /tikz/column 2/.append style={anchor=base east},
    /tikz/column 3/.append style={anchor=base west} ]
    \varphi : & E \arrow[r]  & E'\,, \\
    & P \arrow[r, mapsto] & (x'(P), y'(P))\,.
  \end{tikzcd}
\end{equation*}
 the quotient by $\gT$
isogeny. This is a degree two separable isogeny.
The quotient curve $E'$ has  Weierstrass equation
\begin{displaymath}
y'^2+a_1x'y'+a_3y'=x'^3+a_2x'^2+A_4x'+A_6.
\end{displaymath}
 The meaning of the lemma below is that the local parameters
$z$ and $z'$ are very close in the neighborhood of  the origin.
\begin{lemma}\label{lemma:z'}
  In the neighborhood  of $O$ we have $z'=z+O(z^4)$
  and $\frac{1}{z'}=\frac{1}{z}+O(z^2)$.
\end{lemma}
\begin{myproof}
  According to Vélu's formulae (\ref{eq:vel1})
  \[z'=-\frac{x'}{y'}=z\frac{1+\frac{x_T-x(T)}{x}}{1+\frac{y_T-y(T)}{y}}
  =z(1+O(z^3)).\]
\end{myproof}

\section{Evaluation using elliptic butterflies}\label{sec:but}

Under the hypotheses
and  notation
of Sections~\ref{sec:ebfu} and \ref{sec:2},
we let   $b$ be
a point in $E(\KK)$ such
that \[db\not = O.\]
We set \[b'=\varphi(b)\text{\,\,\, and \,\,\,} d'=2^{\delta-1}=d/2\] and we
notice that \[\hat\varphi (d'b')\not =O\]
\smallskip
where $\hat\varphi$ is the dual isogeny to $\varphi$.
Given $f$ in $\cL (\gt)$ by its coordinates
in  basis $u$, we want to evaluate $f$
at every point in the coset $b\,  + \gt$.
We reduce this problem to two similar problems
on $E'$ each of halved size. We show that
this reduction is achieved at the expense
of $\cQ d$ additions and multiplications by a constant in $\KK$.
Using this recursion, the evaluation of $f$
at $b\,  + \gt$ is achieved at the expense
of $\cQ d\log (d)$ additions and multiplications by a constant in $\KK$.
The reduction step is the elliptic analogue
of the butterfly diagram appearing in the
standard FFT.

We call
  $O'=(0,1,0)$ the origin on $E'$.
We denote by  \[t' =\varphi(t)\] a point of order $d'$ on $E'$.
The group $\gtp$
can be seen as an effective  divisor on $E'$ with degree $d'$.
We consider  the linear space $\cL(\gtp)$  associated with $\gtp$.
In case $d'\not = 1$ we denote by \[u'=(u'_l)_{l\in \ZZ/d'\ZZ}\text{\,\,\,  and \,\,\,} v'=(v'_l)_{l\in \ZZ/d'\ZZ}\]
the two bases of $\cL (\gtp)$ as constructed   in Section \ref{sec:ebfu}.
For $l$ in $\ZZ/d'\ZZ$
\[u'_l=u_{lt',(l+1)t'}+(1-\agot')/d'\]
where
\[\agot'= -a_1+\Gamma(O,t',2t')+\Gamma(O,2t',3t')+\Gamma(O,3t',4t')+\dots
+\Gamma(O,(d'-2)t',(d'-1)t').\]
In case $d'=1$ we set $u'_0=v'_0=1$ by convention.
We say that a function $f$ on $E$ is \[T\text{-even if }f\circ \tau_T=f\text{\,\,\,
and \,\,\,}T\text{-odd  if }f\circ \tau_T=-f.\] We denote by
\begin{equation}\label{eq:f+-}
  f^+=\frac{f+f\circ \tau_T}{2}\text{\,\,\, and \,\,\,}f^-=\frac{f-f\circ \tau_T}{2}\end{equation}
the $T$-even and $T$-odd parts of $f$.
We shall need a small degree $T$-odd function on $E$. We let
\[\theta = u_{0,T}+\frac{a_1}{2}\] and check that
$\theta$ is $T$-odd. We deduce that the divisor of
$\theta$ is
\[(\theta)=[U]+[U+T]-[O]-[T]\]
where $U$ is any $2$-torsion point on $E$ not in $\gT$.
As a consequence, the product of $\theta$
and its translate $\theta\circ \tau_U$
is a non-zero constant.
We write \[f^0=f^+\text{\,\,\, and \,\,\,}f^1=\theta  f^-\] and we
check that \[f=f^0+\theta^{-1}f^1\]
where both $f^0$ and $f^1$  are $T$-even and thus can be seen
as functions on $E'$.
Indeed
\[f^0 \in \cL (\gtp) \text{\,\,\, and \,\,\,}
f^1 \in \cL (\gtp +O'-U')\]  where
\[U'=\varphi(U).\]
Assume that $f$ is given in the basis $u=(u_l)_{l\in \ZZ/d\ZZ}$ so
  \[f=\sum_{l\in \ZZ/d\ZZ}f_l\, u_l\]
  and
  \begin{equation*}
    f^0=f^+=\frac{1}{2}\sum_{l\in \ZZ/d\ZZ }(f_l+f_{l+d'})u_l
  =\frac{1}{2}\sum_{0\leqslant l \leqslant d'-1}(f_l+f_{l+d'})(u_l+u_{l+d'}).\end{equation*}
From Equation (\ref{eq:TuAB}) and Lemma \ref{lemma:z'} we deduce   that
there is a constant $\ugot$
such that for every integer
  $l$ in $[0, d'-1]$
  \begin{equation}\label{eq:tgot2}
    {u_l+u_{l+d'}}=u'_l+\ugot\,.\end{equation}
    Summing out over $l$ we see that $\ugot=0$. So
            \begin{equation}\label{eq:f+}
          f^0=f^+=\frac{1}{2}\sum_{0\leqslant l \leqslant
            d'-1}(f_l+f_{l+d'})u'_l\,,
        \end{equation}
and    we can compute the coordinates of $f^0$
in the basis $u'$ at the expense of $\cQ d$ additions and multiplications by a constant
in $\KK$.
Similarly
  \begin{equation}\label{eq:f1th}f^1=\theta  f^-=\frac{1}{2}\sum_{l\in \ZZ/d\ZZ }(f_l-f_{l+d'})\cdot \theta \cdot  u_l
    =\frac{1}{2}\sum_{0\leqslant l \leqslant d'-1}(f_l-f_{l+d'})\cdot \theta
    \cdot (u_l-u_{l+d'}).\end{equation}
For every integer
$l$ in $[1, d'-2]$ there exist constants
$\bgot_l$ and $\cgot_{l+1}$ in $\KK$ such that
\begin{equation}\label{eq:chb}
  {\theta \cdot  (u_l-u_{l+d'})}=\bgot_l\cdot (v'_l-v'_l(U'))
  +\cgot_{l+1}\cdot (v'_{l+1}-v'_{l+1}(U')).\end{equation}
      In the special case $l=0$ we find  constants
      $\bgot_\ast$ and $\cgot_1$, such that
      \begin{equation}\label{eq:chb2}
        {\theta  \cdot (u_0-u_{d'})}=  \bgot_\ast \cdot
        (x'-x'(U'))  + \cgot_1\cdot (v'_1-v'_1(U')).\end{equation}
      In the special case $l=d'-1$ we find two
      constants    $\bgot_{d'-1}$ and $\cgot_\ast$, such that
      \begin{equation}\label{eq:chb3}
        \theta \cdot (u_{d'-1}-u_{d-1})=\bgot_{d'-1} \cdot (v'_{d'-1}-v'_{d'-1}(U'))+
        \cgot_\ast \cdot (x'-x'(U')).\end{equation}
When $d'=1$, Equations (\ref{eq:chb}),
      (\ref{eq:chb2}) and (\ref{eq:chb3}) simplify to
  \begin{equation}
    \label{eq:deq2chb}
    \theta \cdot (u_{0}-u_{1}) =  (\bgot_\ast + \cgot_\ast)(x'-x'(U'))\,.
  \end{equation}
      We deduce that $f^1=f^2+f^3$ where
      \begin{align}\label{eq:f2}
        \begin{split}
          f^2&=\frac{1}{2}\sum_{1\leqslant l \leqslant d'-1}(f_l-f_{l+d'})\cdot
          \bgot_l \cdot (v'_l-v'_l(U'))\\[1em]&+  \frac{1}{2}\sum_{0\leqslant l \leqslant d'-2}(f_l-f_{l+d'})\cdot \cgot_{l+1} \cdot (v'_{l+1}-v'_{l+1}(U'))\end{split}\end{align}
      and \begin{equation}\label{eq:f3}f^3=\frac{1}{2}\left( (f_0-f_{d'}) \bgot_\ast
        + (f_{d'-1}-f_{d-1})\cgot_\ast\right) (x'-x'(U')).\end{equation}

      The function $f^2$ belongs to $\cL(\gtp-U')$ and  according to Equation (\ref{eq:f2}) one can compute its  coordinates in the basis $v'$ at the expense
      of  $\cQ d$ additions and multiplications by a constant in $\KK$.
      As explained in Section \ref{sec:ebfu} we can deduce  the coordinates of $f^2$ in the basis $u'$ at the expense of  $\cQ d$ more
such      operations in $\KK$.
When $d'\not = 1$ we  evaluate $f^0$ and $f^2$  recursively as functions on $E'$.
In the special case when $d'=1$ these two   functions
are given constants.
      Using Equation (\ref{eq:f3}) we evaluate $f^3$ at the expense of
      $3$ additions and $d'+2$ multiplications
      provided we have precomputed the values of
      $x'-x'(U')$ at $b'+\gtp$.  We deduce the values of $f^1=f^2+f^3$ at the expense of $d'$ additions.
      We deduce the values of $f=f^0+\theta^{-1}  f^1$ at the expense of $d$ multiplications and $d$ additions.
     As for the constants appearing in Equations (\ref{eq:f2}) and (\ref{eq:f3}) we deduce
      from the examination of Taylor expansions
      that
      \begin{align}\label{eq:bcgot}
        \begin{split}
          \bgot_\ast &=\cgot_\ast=1,\\[1em]
          \bgot_l&=-\Gamma(O,T,lt)-a_1/2 =-\theta (lt)\text{\,\,\, for \,\,\,}   l \in [1, d'-1],\\[1em]
          \cgot_l&=\Gamma(O,T,lt)+a_1/2 =\theta (lt)\text{\,\,\, for \,\,\,}   l \in [1, d'-1].
        \end{split}
      \end{align}

\bigskip

      \begin{figure}[htbp]
        \begin{linebox}[width=0.95\linewidth]{Algorithm}
          \begin{small}
            \begin{algorithmic}
              \Function{Butterfly\_Evaluate\,}{$\vec{f}$}
              \State \hspace{0.5cm}$\triangleright$ Constants:
              \State \hspace{1cm}$\_$ The constants $\vec{\agot}\,'$,
              $\vec{\bgot}$ and $\vec{\cgot}$ are defined by
              $\bgot_0=\bgot_\ast$\,, $\cgot_0=\cgot_\ast$ and by
              Eq.~\eqref{eq:agot} and~\eqref{eq:bcgot}\,.
              \State \hspace{1cm}$\_$ The constant $\vec{\xgot}\,'$ denotes
              $(x'(b'+l\,t') - x'(U'))_{l \in \ZZ/d'\ZZ}$\,.
              \State \hspace{1cm}$\_$ The constant $\vec{\vgot}\,'$ denotes
              $(v'_l(U'))_{l \in \ZZ/d'\ZZ}$\,.
              \State \hspace{1cm}$\_$  The constant $\vec{\tgot}$ denotes $
              (\theta(b+l\,t))_{l=0,\,\ldots,\, d'-1}$\,.
              \State \hspace{1cm}$\_$ The constant $\vec{\mgot}\,'$ denotes
              $\vec{\bgot} \,\diamond\, \vec{\vgot}\,'$ and $\vec{\ngot}\,'$
              denotes
              $\vec{\mgot}\,' + \sigma^{-1}( \vec{\cgot} \,\diamond\,
              \vec{\vgot}\,' )$\,.
              \State
              \State \algorithmicif\ $d = 1$ \ \algorithmicthen\
              \Return $\vec{f}$\,; \ \algorithmicend\ \algorithmicif
              \Comment{Recursion end}
              \State $\vec{f}\,^{\text{low}} \gets \vec{f}\,_{[\,0,\ d'-1\,]}$\,; \
              $\vec{f}\,^{\text{high}} \gets \vec{f}\,_{[\,d',\ d-1\,]}$ \Comment{Split $\vec{f}$ in half}
              \State
              $\vec{f}\,^{+} \gets ({1}/{2}) \cdot (\vec{f}\,^{\text{low}} +
              \vec{f}\,^{\text{high}})$\,; \
              $\vec{f}\,^{-} \gets ({1}/{2}) \cdot (\vec{f}\,^{\text{low}} -
              \vec{f}\,^{\text{high}})$ \Comment{Symmetrization}
              \State $\vec{\alpha}\,^{+} \gets \textsc{Butterfly\_Evaluate}\,(\vec{f}\,^{+})$
              \Comment{Recursive call}
              \State $r \gets  \bgot_{0} \, {f_{0}}^{-}  +  \cgot_{0}
              \,{f}_{d'-1}^{-} $\,; $t \gets 0$
              \Comment{Constants}
              \If {$d > 2$}
              \State $t \gets \mgot_{0}' \, {f}_{0}^{-} + (\ngot_{d'-1}' - {\mgot}_{d'-1}')
              \, {f}_{d'-1}^{-} - \langle\vec{\ngot}\,',\ \vec{f}\,^{-}\rangle$
              \EndIf
              \State
              $\vec{f}\,^{-} \gets \textsc{BiDiagonal\_Evaluate}\,(\vec{\bgot},\, \vec{\cgot},\,
              \vec{f}\,^{-})$\,;\ ${f}_{0}^{-} \gets t$
              \Comment{Bi-diagonal formulas}
              \State $\vec{f}\,^{-} \gets \textsc{VtoU\_BaseChange}\,(\vec{\agot}\,',\, \vec{f}\,^{-})$
              \Comment{Base change from $(v'_l)$ to $(u'_l)$}
              \State $\vec{\alpha}\,^{-} \gets \textsc{Butterfly\_Evaluate}\,(\vec{f}\,^{-})$
              \Comment{Recursive call}
              \State
              $\vec{\alpha}\,^{-} \gets (\, \vec{\alpha}\,^{-} + r \cdot \,\vec{\xgot}\,'\,)\, \diamond\, \vec{\tgot}\,^{\mathrm{inv}}$
              \Comment{$\theta$-multiplication}
              \State \Return $( \vec{\alpha}\,^{+} + \vec{\alpha}\,^{-} ) \oplus (
              \vec{\alpha}\,^{+} - \vec{\alpha}\,^{-} )$
              \EndFunction
            \end{algorithmic}
          \end{small}
        \end{linebox}
        \caption{Butterfly Evaluation}
        \label{alg:butterfly_evaluate}
      \end{figure}

      \begin{proposition}[Fast elliptic evaluation]\label{prop:eva}
  There exists a constant $\cQ$ such that the following is true.
        Let $\KK$ be a field with odd characteristic.
        Let $E$ be an elliptic curve over $\KK$.
        Let $\delta \geqslant 1$ be an integer.
        Let $t$ be a point of order $d=2^\delta$
        in $E(\KK)$.
        Let $b$ be a point in $E(\KK)$ such that
        $db\not=O$.
        There exists a straight-line program that on input
a function  $f$ in $\cL (\gt)$  given by its coordinates
        in either basis $u$ or $v$ defined in Section
        \ref{sec:ebfu}, computes the values $f(b+lt)$ for
        $l\in [0, d-1]$ at the expense of $\cQ d\log(d)$
        operations in $\KK$.
      \end{proposition}
\begin{myproof}      This results from the discussion above.
      For clarity, the evaluation procedure is presented as concise pseudo-code in
      Figure~\ref{alg:butterfly_evaluate}.\end{myproof}

\section{Interpolation}\label{sec:int}

In the context of the beginning of Section \ref{sec:but}
and using a similar recursion we explain how to recover
the coefficients of a function $f\in \cL (\gt)$
in the basis $u$ from its values at  $b\,  + \gt$.
We assume that we are given $d$ scalars $(\alpha_l)_{0\leqslant
  l\leqslant d-1}$
in $\KK$.
We want to compute the coordinates $(f_l)_{l\in \ZZ/d\ZZ}$ in
the basis $u$ of the unique function $f$ in $\cL (\gt)$
such that \[f(b+lt)=\alpha_l \text{\,\, for every \,\,} l\in  [0,d-1].\]
We use a recursion again.
We first compute \[\alpha_l^+=\frac{\alpha_l+\alpha_{l+d'}}{2}
\text{ \,\,\, for \,\,\,} 0\leqslant l\leqslant d'-1\]
the values of $f^+$ at $b'\,  + \gtp$.
If $d'\not = 1$ we deduce, by recursion,
the coordinates $(f^+_l)_{0\leqslant l\leqslant d'-1}$
of $f^+$ in the basis $u'$.
In the special case when $d'=1$ we simply
have $f^+_0=\alpha^+_0$.
From Equations~(\ref{eq:f+}) we deduce
      \begin{align*}
        \begin{split}
          f^+&=\sum_{0\leqslant l\leqslant d'-1}f_l^+ \cdot u'_l\\[1em]
          &=\frac{1}{2}\sum_{0\leqslant l\leqslant d'-1}(f_l+f_{l+d'})u'_l.\end{split}\end{align*}
      So for $0\leqslant l\leqslant d'-1$
      \begin{equation*}\label{eq:f+pap}
f_l^+=\frac{f_l+f_{l+d'}}{2}.
        \end{equation*}
     This gives half the information we need to quickly compute the
      $(f_l)_{l\in \ZZ/d\ZZ}$. We will be done when   we compute
the
      \[\frac{f_l-f_{l+d'}}{2}\text{\,\,\, for\,\,\,}  0\leqslant l\leqslant d'-1.\]
      According to Equation (\ref{eq:f1th}) these $d'$ scalars are the coordinates of $f^1=\theta \cdot f^-$
      in the basis of $\cL(\gtp +O'-U')$ made of the
      $\theta \cdot (u_l-u_{l+d'})$ for $l\in [0,d'-1]$.
We will use an auxiliary function
      $\xi_b$ on $E'$ having degree $d'+1$,
      polar divisor $\gtp+O'$, and vanishing
      at $b'+lt'$ for every $l\in [0,d'-1]$. Such a function is unique
      up to  a multiplicative constant, and it has no pole
      and no zero at $U'$ because $db\not =O$. We  therefore assume that
      $\xi_b (U')=1$. There exists
      $d'+1$ scalars $\dgot_0$, $\dgot_1$, \ldots, $\dgot_{d'-1}$,
      $\dgot_\ast$ such that \begin{equation}\label{eq:xi}
        \xi_b=\dgot_\ast \cdot x'+\sum_{0\leqslant l \leqslant d'-1}\dgot_l\cdot v'_l=
        1+\dgot_\ast \cdot (x'-x'(U'))+ \sum_{1\leqslant l \leqslant d'-1}\dgot_l \cdot (v'_l -v'_l(U')).\end{equation}
      We assume that we have precomputed these $d'+1$ scalars.
         We now compute \[\alpha_l^-=\frac{\alpha_l-\alpha_{l+d'}}{2}
     \text{\,\,\, and \,\,\,} \alpha_l^1=\theta(b+lt) \cdot \alpha_l^-
     \text{ \,\,\, for \,\,\,} 0\leqslant l\leqslant d'-1\]
     at the expense of $\cQ d$ operations in $\KK$,
     assuming we have precomputed the values of $\theta$
     at $b+\gt$. The $\alpha_l^1$
     are the values of $f^1$ at  $b'+\gtp$.
     We let
 $f^\star$ be the       function  in $\cL(\gtp)$
     that takes value $\alpha_l^1$ at $b'+lt$.
     We obtain by recursion the coordinates
     of $f^\star$ in the basis $u'$. We base change to $v'$ and obtain
          $d'$ scalars $f^\star_l$
     such that \[f^\star = \sum_{0\leqslant l \leqslant d'-1}f^\star_l v'_l.\]
      We are not quite
     done. The function $f^\star$ is not the one we are looking
     for because it does not vanish at $U'$. So we
     compute \[f^\star(U')=\sum_{0\leqslant l \leqslant d'-1}f^\star_l\cdot v'_l(U')\]
     at the expense of $d'$ multiplications
     and $d'-1$ additions, provided we
     have precomputed the $v'_l(U')$. And we
     write
     \[f^\star = f^\star(U')
    +\sum_{1\leqslant l \leqslant d'-1}f^\star_l\cdot (v'_l-v'_l(U'))\,.\]
    We deduce
      \begin{align*}
        \begin{split}
          f^1&=f^\star-f^\star(U')\cdot \xi_b\\[1em]
          &=-f^\star(U')\cdot \dgot_\ast \cdot (x'-x'(U'))+
          \sum_{1\leqslant l \leqslant d'-1}\big(f^\star_l-\dgot_l\cdot f^\star(U')\big)\cdot (v'_l-v'_l(U')).
        \end{split}
      \end{align*}
      This gives the coordinates
      \[(s_\ast, s_1, s_2, \ldots, s_{d'-1})\]
      of $f^1$ in the basis
      of $\cL(\gtp +O'-U')$ made of $x'-x'(U')$ and
      the $v'_l-v'_l(U')$ for $l\in [1,d'-1]$.
      Indeed \[s_l = f^\star_l-\dgot_l\cdot f^\star(U')\text{\,\,\, for \,\,\,} 1\leqslant l\leqslant d'-1, \text{\,\,\, and \,\,\,} s_\ast = -f^\star(U')\cdot \dgot_\ast.\]
      We want to  deduce
      the coordinates
      \[(t_0, t_1, t_2, \ldots, t_{d-1})\] of $f^1$
      in the  basis of
      $\cL(\gtp +O'-U')$ made of the
      $\theta \cdot (u_l-u_{l+d'})$ for $l\in [0,d'-1]$.
According
      to Equations (\ref{eq:chb}),
      (\ref{eq:chb2}), and (\ref{eq:chb3}) this amounts
      to solving the following linear system.

      \[\begin{pmatrix}
            \bgot_\ast & 0 & \cdots & 0 & \cgot_\ast\\
      \cgot_1 & \bgot_1 & 0 & \cdots & 0\\
      0  & \cgot_2 & \bgot_2 & \ddots & \vdots\\
      \vdots & \ddots & \ddots & \ddots & 0 \\
      0 &\cdots & 0 & \cgot_{d'-1} & \bgot_{d'-1}
      \end{pmatrix}
      \begin{pmatrix}t_0\\t_1\\\vdots \\t_{d-2}\\t_{d'-1}\end{pmatrix}
        =
              \begin{pmatrix}s_\ast\\s_1\\\vdots \\s_{d'-2}\\s_{d'-1}\end{pmatrix}\]

                \smallskip
                Solving this system requires $\cQ d$ operations in $\KK$ as
      recalled in section \ref{sec:bidia}.
      In the special case when $d'=1$ we simply invoke Equation
      (\ref{eq:deq2chb}).

            \begin{proposition}[Fast elliptic interpolation]\label{prop:int}
  There exists a constant $\cQ$ such that the following is true.
        Let $\KK$ be a field with odd characteristic.
        Let $E$ be an elliptic curve over $\KK$.
        Let $\delta \geqslant 1$ be an integer.
        Let $t$ be a point of order $d=2^\delta$
        in $E(\KK)$.
        Let $b$ be a point in $E(\KK)$ such that
        $db\not=O$.
        There exists a straight-line program that on input
        $d$ scalars $(\alpha_l)_{0\leqslant l\leqslant d-1}$,
        computes
        a function  $f$ in $\cL (\gt)$,  given by its coordinates
        in either basis $u$ or $v$ defined in Section
        \ref{sec:ebfu}, such that  $f(b+lt)=\alpha_l$ for every
        $l\in [0, d-1]$,  at the expense of $\cQ d\log(d)$
        operations in $\KK$.
      \end{proposition}
\begin{myproof}      This results from the discussion above.
      To aid understanding, the interpolation procedure is shown in compact pseudo-code in Figure~\ref{alg:butterfly_interpolate}.
      \end{myproof}

\begin{figure}[htbp]
  \begin{linebox}[width=0.95\linewidth]{Algorithm}
    \begin{small}
      \begin{algorithmic}
        \Function{Butterfly\_Interpolate\,}{$\vec{\alpha}$}
        \State \hspace{0.5cm}$\triangleright$ Constants:
        \State \hspace{1cm}$\_$ The constants $\vec{\agot}\,'$, $\vec{\bgot}$
        and $\vec{\cgot}$ are defined by $\bgot_0=\bgot_\ast$\,,
        $\cgot_0=\cgot_\ast$ and by Eq.~\eqref{eq:agot}
        and~\eqref{eq:bcgot}\,.
        \State \hspace{1cm}$\_$ The constant $\vec{\dgot}$ is defined by
        Eq.~\eqref{eq:xi}.
        \State \hspace{1cm}$\_$ The constant $\vec{\vgot}\,'$ denotes
        $(v'_l(U'))_{l \in \ZZ/d'\ZZ}$\,.
        \State \hspace{1cm}$\_$ The constant $\vec{\tgot}$ denotes $
        (\theta(b+l\,t))_{l=0,\,\ldots,\, d'-1}$\,.
        \State
        \State \algorithmicif\ $d = 1$ \ \algorithmicthen\ \Return
        $\vec{\alpha}$\,; \ \algorithmicend\ \algorithmicif \Comment{Recursion
          end}
        \State
        $\vec{\alpha}\,^{\text{low}} \gets \vec{\alpha}\,_{[\,0,\ d'-1\,]}$\,;
        \ $\vec{\alpha}\,^{\text{high}} \gets \vec{\alpha}\,_{[\,d',\ d-1\,]}$
        \Comment{Split $\vec{\alpha}$ in half}
        \State
        $\vec{\alpha}\,^{+} \gets ({1}/{2}) \cdot (\vec{\alpha}\,^{\text{low}}
        + \vec{\alpha}\,^{\text{high}})$\,; \
        $\vec{\alpha}\,^{-} \gets ({1}/{2}) \cdot (\vec{\alpha}\,^{\text{low}}
        - \vec{\alpha}\,^{\text{high}})$ \Comment{Symmetrization}
        \State
        $\vec{f}\,^{+} \gets
        \textsc{Butterfly\_Interpolate}\,(\vec{\alpha}\,^{+})$
        \Comment{Recursive call}
        \State
        $\vec{f}\,^{-} \gets \textsc{Butterfly\_Interpolate}\,(\vec{\tgot}
        \,\diamond\, \vec{\alpha}\,^{-})$ \Comment{Recursive call}
        \State
        $\vec{f}\,^{-} \gets \textsc{UtoV\_BaseChange}\,(\vec{\agot}\,',
        \vec{f}\,^{-})$ \Comment{Base change from $(u'_l)$ to $(v'_l)$}
        \State
        $f_\ast \gets \langle \vec{\vgot}\,'\,,\,\vec{f}\,^{-}\rangle$\,;%
        \
        $\vec{f}\,^{-} \gets \vec{f}\,^{-} - f_\ast \cdot \vec{\dgot}\,;$%
        \ ${f}_{0}^{-}\gets -{\dgot}_{0}\,f_\ast$
        \Comment{$\zeta_b$-normalization}
        \State
        $\vec{f}\,^{-} \gets \textsc{BiDiagonal\_Invert}(\vec{\bgot},\,
        \vec{\cgot},\, \vec{f}\,^{-})$ \Comment{Bi-diagonal inversion}
        \State \Return
        $( \vec{f}\,^{+} + \vec{f}\,^{-} ) \oplus ( \vec{f}\,^{+} -
        \vec{f}\,^{-} )$
        \EndFunction
      \end{algorithmic}
    \end{small}
  \end{linebox}
  \caption{Butterfly Interpolation}
  \label{alg:butterfly_interpolate}
\end{figure}

                \section{Reduction}\label{sec:redu}
                                Under the hypotheses at the beginning of Section~\ref{sec:ebfu}
and assuming that the base field $\KK$
has odd characteristic, we let   $b$ be
a point on $E$ defined over a  separable closure $\Ks$ of $\KK$. We assume that  $db\not =O$.
We call $B$ the coset of $b$ under the action of $\gt$.
So
\[B=b+\gt=\{b, b+t, b+2t, \ldots, b+(d-1)t\}.\] We assume that $B$
is left invariant
by the absolute Galois  group of $\KK$. This is equivalent
to the image of $b$ by the isogeny with kernel $\gt$
  being $\KK$-rational.
For every $l$ in $\ZZ/d\ZZ$ we denote by $x_l$ the function
defined by \begin{equation}\label{eq:xl}
  x_l=x_{lt}=x\circ \tau_{-lt}.\end{equation}
We assume that  we are given scalars $(F_l)_{l\in \ZZ/d\ZZ}$
in $\KK$ and we consider  the  function
\[F = \sum_{l\in \ZZ/d\ZZ} F_l\,  x_{l}.\]
We are interested in such  functions
because they form a $t$-invariant
supplementary subspace to $\cL (\gt)$
in $\cL (2\gt)$. Given such a function $F$ there exists a unique function
\[f = \sum_{l\in \ZZ/d\ZZ} f_l \, u_{l}\]
such that $f$ and $F$ agree on $B$ that is \[f(b+lt)=F(b+lt) \text{\,\, for every \,\,} l \in \ZZ/d\ZZ.\] We write this condition $f \equiv F \pmod B$.
The purpose of this section is to compute the $(f_l)_{l\in \ZZ/d\ZZ}$
once given the $(F_l)_{l\in \ZZ/d\ZZ}$.
We shall achieve this  at the
expense of $\cQ d\log (d)$ additions and scalar multiplications
thanks to a  recursion as in the previous sections.
We define $f^+$ and $f^-$ as in Equation (\ref{eq:f+-})
and similarly
\begin{equation*}\label{eq:F+-}
  F^+=\frac{F+F\circ \tau_T}{2}\text{\,\,\, and \,\,\,}F^-=\frac{F-F\circ \tau_T}{2}.\end{equation*} We set
\[B'=\varphi (B)=b'+\gtp=\{b', b'+t', b'+2t', \ldots, b'+(d'-1)t'\}\]
and notice that $F^+\equiv f^+\pmod{B}$
and $F^-\equiv f^-\pmod{B}$. From Equation (\ref{eq:vel1})
we deduce
\begin{align}\label{eq:4l}
  \begin{split}
    F^+&=\sum_{l\in \ZZ/d\ZZ}\frac{F_l+F_{l+d'}}{2}\, x_l\\[1em]
    &=
    \sum_{0\leqslant l \leqslant d'-1}\frac{F_l+F_{l+d'}}{2}\, (x_l+x_{l+d'})\\[1em]
    &=
\sum_{0\leqslant l \leqslant  d'-1}\frac{F_l+F_{l+d'}}{2}\, (x'_l+x(T))\\[1em]&=
F^+_\ast+\sum_{0\leqslant l \leqslant d'-1}F^+_l\cdot x'_l \end{split}\end{align}
where \[F^+_\ast=(x(T)/2)\sum_{0\leqslant l \leqslant d-1} F_l\text{\,\,\, and \,\,\,}
F^+_l=\frac{F_l+F_{l+d'}}{2} \text{\,\, for \,\,} 0\leqslant l \leqslant d'-1.\]
When $d'\not = 1$, we
obtain by recursion $d'$ scalars $(f^+_l)_{0\leqslant l \leqslant d'-1}$
such that \[\sum_{0\leqslant l \leqslant d'-1}F^+_l\cdot x'_l
\equiv \sum_{0\leqslant l \leqslant d'-1}f^+_l\cdot u'_l \bmod B'.\]
Using Equation (\ref{eq:tgot2}) we deduce
\begin{align*}
  \begin{split}
F^+
&\equiv  F^+_\ast+\sum_{0\leqslant l \leqslant d'-1}f^+_l\cdot u'_l \bmod B\\[1em]
&=  F^+_\ast   +\sum_{0\leqslant l \leqslant  d'-1}f^+_l\cdot (u_l+u_{l+d'})\\[1em]
&=F^+_\ast \sum_{0\leqslant l \leqslant d-1}u_l+\sum_{0\leqslant l \leqslant d'-1}f^+_l\cdot
(u_l+u_{l+d'})\\[1em]
&= \sum_{0\leqslant l \leqslant d'-1}(f^+_l+F^+_\ast)(u_l+u_{l+d'}).
    \end{split}\end{align*}
So we have reduced $F^+$ in that case. In the special case when  $d'=1$ we notice that
\begin{equation}\label{eq:d'=1}
  x'_0=x'\equiv x'(b')\pmod {B'}\end{equation} and we deduce from Equation
(\ref{eq:4l}) that
\[F^+\equiv \big(F_\ast^++F_0^+\cdot x'(b')\big)\cdot (u_0+u_1)\pmod B.\]
We now
proceed to reducing $F^-$. We define
\begin{align}\label{eq:F1t}
  \begin{split}
    F^1&=\theta \cdot F^-\\[1em]
    &=\sum_{l\in \ZZ/d\ZZ}\frac{F_l-F_{l+d'}}{2}\cdot \theta \cdot x_l\\[1em]
    &=\sum_{0\leqslant l \leqslant d'-1}\frac{F_l-F_{l+d'}}{2}\cdot \theta \cdot (x_l-x_{l+d'}).
\end{split}\end{align}
For every integer
$l$ in $[1, d'-1]$ the product $\theta \cdot  (x_l-x_{l+d'})$
is a function on $E'$. It belongs to the linear space
$\cL(-2[lt']-[O']+[U'])$. So   there exist constants
$\egot_l$ and $\fgot_l$ in $\KK$ such that
\begin{equation}\label{eq:chbx}
  {\theta \cdot  (x_l-x_{l+d'})}=\fgot_l \cdot (x'_l-x'_l(U'))
  +\egot_l \cdot (v'_{l}-v'_{l}(U')).\end{equation}
In the special case $l=0$ the product
 $\theta \cdot  (x_0-x_{d'})$
belongs to the linear space
$\cL(-3[O']+[U'])$. So
there exist   constants
      $\fgot_0$ and $\ggot$, such that
      \begin{equation}\label{eq:chbx2}
        {\theta \cdot (x_0-x_{d'})}=  \fgot_0\cdot
(x'-x'(U'))  + \ggot \cdot (y'-y'(U')).\end{equation}
      Using Equations (\ref{eq:F1t}), (\ref{eq:chbx}), and (\ref{eq:chbx2}) we write
\begin{align}\label{eq:F1sp}
      \begin{split}
    F^1&=F^2+F^3+F^4 \text{\,\,\, with\,\,\,}\\[1em]
    F^2&=\frac{1}{2}\sum_{1\leqslant l \leqslant d'-1}(F_l-F_{l+d'})\cdot
    \egot_l \cdot v'_l  \\[1em]&-\frac{1}{2}\sum_{1\leqslant l \leqslant d'-1}(F_l-F_{l+d'})\cdot \big(\fgot_l\cdot x'_l(U')+\egot_l\cdot v'_l(U')\big)\\[1em]
       &-\frac{1}{2}(F_0-F_{d'})\cdot \big(\fgot_0\cdot x'_0(U')+\ggot \cdot y'(U')\big),\\[1em]
    F^3&=\frac{1}{2}\sum_{0\leqslant l \leqslant d'-1}(F_l-F_{l+d'})\cdot \fgot_l\cdot x'_l,\\[1em]
    F^4&=\frac{1}{2}(F_0-F_{d'})\cdot  \ggot \cdot y'.
\end{split}\end{align}
The above
expression for $F^2$
provides its coordinates $(F^2_l)_{0\leqslant l\leqslant d'-1}$
 in the basis $v'$ of $\cL(\gtp)$.
There exist
scalars $(F^3_l)_{0\leqslant l\leqslant d'-1}$ such that
\[F^3\equiv \sum_{0\leqslant l\leqslant d'-1}F^3_l\cdot v'_l \pmod{B'}.\]
To  obtain these scalars
we  use the
expression for $F^3$ in Equation (\ref{eq:F1sp}).

\noindent If $d'\not = 1$ we
apply reduction recursively, then base change from $u'$ to $v'$.

\noindent In the special case when $d'=1$
we deduce from Equation (\ref{eq:d'=1})  that \[F^3\equiv \frac{1}{2}(F_0-F_1)
\cdot \fgot_0\cdot x'(b') \cdot v'_0\pmod {B'}.\]
As for $F^4$, we assume that we have precomputed scalars
$(\hgot_l)_{0\leqslant l\leqslant d'-1}$ such that
\begin{equation}\label{eq:hgot}
  y'\equiv \sum_{0\leqslant l\leqslant d'-1}\hgot_l\cdot v'_l \pmod{B'}.
\end{equation}
We compute \[F^4_l=(1/2)(F_0-F_{d'})\cdot \ggot \cdot \hgot_l\text{\,\, for \,\,}
l\in [0,d'-1]\]
and we have
\[F^4\equiv \sum_{0\leqslant l\leqslant d'-1}F^4_l\cdot v'_l \pmod{B'}.\]
Finally  we compute \[f^\star_l=F^2_l+F^3_l+F^4_l.\]The reduction
of $F^1$ modulo $B'$ is \[f^\star =
\sum_{0\leqslant l\leqslant d'-1}f^\star_l\cdot v'_l \equiv F^1 \pmod{B'}.\]
Since $F^1=\theta \cdot  F^-$ we would now  divide $f^\star$  by $\theta$. But $\theta$ vanishes at
$U$ and $f^\star$ does not. Again we overcome this difficulty
thanks to  the function
$\xi_b$ introduced before Equation (\ref{eq:xi}).
We compute
\[f^\star (U')=\sum_{0\leqslant l\leqslant d'-1}f^\star_l\cdot v'_l(U')\]
and check that the function
\[f^1 = f^\star - f^\star (U')\cdot \xi_b \in \cL (\gtp +O'-U')\]
is congruent to $F^1$ modulo $B'$. So
\begin{align}\label{eq:F-t}
      \begin{split}
        f^-&=f^1\cdot \theta^{-1}\in \cL(\gt)\\[1em]
           &=f^\star\cdot \theta^{-1} - f^\star (U')\cdot \xi_b\cdot \theta^{-1}\\[1em]
        &=f^\star (U')\cdot (1-\xi_b)\cdot \theta^{-1}+\sum_{1\leqslant l\leqslant d'-1}f^\star_l\cdot \big(v'_l-v'_l(U')\big)\cdot
        \theta^{-1} \in \cL(\gt).
      \end{split}\end{align}
Indeed this function is congruent to $F^-$ and it belongs to $\cL(\gt)$.
We are not quite done because $f^-$
is given as a linear combination of
    $(1-\xi_b)\cdot \theta^{-1}$ and the $(v'_l-v'_l(U'))\cdot \theta^{-1}$.
Fortunately, for every $l$ in $[1,d'-1]$ there exist
    constants $\igot_l$
    and $\jgot_l$ such that
    \begin{equation}\label{eq:passv}
      (v'_l-v'_l(U'))\cdot \theta^{-1}= \igot_l \cdot (v_l-v_{l+d'})+\jgot_l.
    \end{equation}
    And there exist constants $(\lgot_l)_{0\leqslant l\leqslant d-1}$
    such that
    \begin{equation}\label{eq:passxi}
(1-\xi_b)\cdot \theta^{-1}=\sum_{0\leqslant l\leqslant d-1}\lgot_l\cdot v_l.
\end{equation}
Substituting Equations (\ref{eq:passv}) and (\ref{eq:passxi})
in (\ref{eq:F-t}) we terminate the reduction of $F^-$.
We finally compute $f=f^++f^-$.

To complete  this discussion
we give simple expressions for  the constants
involved in Equations (\ref{eq:chbx}),  (\ref{eq:chbx2}), (\ref{eq:passv}).
Thus $\ggot=1$, $\fgot_0=a_1/2$,  and for $l \in [1, d'-1]$
      \begin{align}\label{eq:efijgot}
        \begin{split}
          \egot_l&= x(lt)-x(lt+T),\\[1em]
          \fgot_l&=\theta (lt),\\[1em]
          \igot_l&= \theta^{-1}(lt),\\[1em]
          \jgot_l&=1.
\end{split}
        \end{align}

\bigskip

            \begin{proposition}[Fast elliptic reduction]\label{prop:red}
  There exists a constant $\cQ$ such that the following is true.
        Let $\KK$ be a field with odd characteristic.
        Let $E$ be an elliptic curve over $\KK$.
        Let $\delta \geqslant 1$ be an integer.
        Let $t$ be a point of order $d=2^\delta$
        in $E(\KK)$.
         Let   $b$ be
         a point on $E$ defined over a  separable extension
         of $\KK$. Assume   $db\not = O$.
Let
\[B=b+\gt=\{b, b+t, b+2t, \ldots, b+(d-1)t\}.\] Assume  $B$
is left invariant
by the absolute Galois  group of $\KK$.
        So $B$ is  a reduced $\KK$-subscheme of $E$
        having dimension $0$ and degree $d$.
        Let $(x_l)_{0\leqslant l\leqslant d-1}$  be the functions
        on $E$ defined in Equation \ref{eq:xl}.
        There exists a straight-line program that on input
        $d$ scalars $(F_l)_{0\leqslant l\leqslant d-1}$,
        computes
        $d$ scalars $(f_l)_{0\leqslant l\leqslant d-1}$,
        such that the functions
        \[F = \sum_{l\in \ZZ/d\ZZ} F_l\,  x_{l}
\text{\,\,\,  and \,\,\,} f= \sum_{l\in \ZZ/d\ZZ} f_l\,  u_{l}\]
are congruent modulo $B$, meaning \[f(b+lt)=F(b+lt)
\text{\,\,    for every integer \,\,}l\in [0, d-1].\]
        This straight-line program consists of  no
        more than  $\cQ d\log(d)$
        operations in $\KK$.
      \end{proposition}
\begin{myproof}      This results from the discussion above.
      To make the procedure clearer, we present the reduction routine as brief pseudo-code in Figure~\ref{alg:butterfly_reduction}.
      \end{myproof}

A consequence of Propositions \ref{prop:red}
and \ref{prop:eva} is that, in  case when $b$ is a $\KK$-rational
point, we can evaluate $F= \sum_{l\in \ZZ/d\ZZ} F_l\,  x_{l}$
at all $b+lt$ by first
computing the reduction $f$ of $F$ modulo $B$, then
evaluating $f$ at the  $b+lt$. This takes time
$\cQ d\log( d)$.

\begin{figure}[htbp]
  \begin{linebox}[width=0.95\linewidth]{Algorithm}
    \begin{small}
      \begin{algorithmic}
        \Function{Butterfly\_Reduce}{$\vec{F}$}
        \State \hspace{0.5cm}$\triangleright$ Constants:
        \State \hspace{1cm}$\_$ The constants $\vec{\agot}\,'$ is defined by
        Equation~\eqref{eq:agot} and $\vec{\agot}\,^1$ denotes
        $\vec{\agot}\,_{[\,d',\ d-1\,]} - \vec{\agot}\,_{[\,0,\ d'-1\,]}$\,.
        \State \hspace{1cm}$\_$ The constants $\vec{\egot}$, $\vec{\fgot}$ and
        $\vec{\igot}$ are defined by Equation~\eqref{eq:efijgot}.
        \State \hspace{1cm}$\_$ The constant $\vec{\hgot}$ is defined by
        Equation~\eqref{eq:hgot}.
        \State \hspace{1cm}$\_$ The constant $\vec{\lgot}$ is defined by
        Equation~\eqref{eq:passxi}, restricted to $l\in [0,\, d'-1]$
        and\newline
        \makebox[\linewidth][r]{$\lgot_\ast$ denotes the constant $\lgot_{d'}$.}
        \State \hspace{1cm}$\_$ The constant $\vec{\vgot}\,'$ denotes
        $(v'_l(U'))_{l \in \ZZ/d'\ZZ}$\,.
        \State \hspace{1cm}$\_$ The constant $\vec{\pgot}$ denotes the vector
        $\vec{\egot} \,\diamond\, \vec{\vgot}\,' + \vec{\fgot} \,\diamond\,
        (x'(U'+l\,t'))_{l \in \ZZ/d'\ZZ} $\,, and\newline
        \makebox[\linewidth][r]{ $\pgot_\ast$ denotes the
          constant $\pgot_0 - \fgot_0\,x'(U') - y'(U')$\,.}
        \State
        \State \algorithmicif\ $d = 1$ \ \algorithmicthen\ \Return
        $x'(b') \cdot \vec{F}$\,; \ \algorithmicend\ \algorithmicif
        \Comment{Recursion end}
        \State $\vec{F}\,^{\text{low}} \gets \vec{F}\,_{[\,0,\ d'-1\,]}$\,; \
        $\vec{F}\,^{\text{high}} \gets \vec{F}\,_{[\,d',\ d-1\,]}$
        \Comment{Split $\vec{F}$ in half}
        \State
        $\vec{F}\,^{+} \gets ({1}/{2}) \cdot (\vec{F}\,^{\text{low}} +
        \vec{F}\,^{\text{high}})$\,; \
        $\vec{F}\,^{-} \gets ({1}/{2}) \cdot (\vec{F}\,^{\text{low}} -
        \vec{F}\,^{\text{high}})$ \Comment{Symmetrization}
        \State ${F}_\ast^{+} \gets x(T) \,\sum_{l = 0}^{d'-1}\,{F}_l^{+}$
        \Comment{Star constant}
        \State
        $\vec{f}\,^{+} \gets {F}_\ast^{+} +
        \textsc{Butterfly\_Reduce}\,(\vec{F}\,^{+})$ \Comment{Recursive call
          and shift by ${F}_\ast^{+}$}
        \State
        $\vec{f}\,^{1} \gets \textsc{Butterfly\_Reduce}\,(\vec{\fgot}
        \,\diamond\, \vec{F}\,^{-})$ \Comment{Recursive call}
        \State
        $\vec{f}\,^1 \gets \textsc{UtoV\_BaseChange}\,(\vec{\agot}\,',
        \vec{f}\,^{1})$ \Comment{Base change from $(u'_l)$ to $(v'_l)$}
        \State
        $\vec{f}\,^{1} \gets \vec{f}\,^{1} + \vec{\egot} \,\diamond\,
        \vec{F}\,^{-} + {F}_{0}^{-} \cdot \vec{\hgot}$
        \State ${f}^{1}_{0} \gets {f}^{1}_{0} + \pgot_\ast\, {F}^{-}_{0} %
        - \langle \vec{\pgot} \,,\, \vec{F}\,^{-}\rangle $
        \Comment{Normalization}
        \State
        $f_\ast^{-} \gets \langle \vec{\vgot}\,' \,,\, \vec{f}\,^{1} \rangle$
        \State
        $\vec{f}\,^{-} \gets -\vec{f}\,^{1} \,\diamond\, \vec{\igot} -
        f_\ast^{-} \cdot \vec{\lgot}\,;\ {f}^{-}_{0} \gets \lgot_{\ast}\,
        f_\ast^{-}$ \Comment{Recover $f^{-}$}
        \State
        $\vec{f}\,^{-} \gets \textsc{VtoU\_BaseChange}\,(\vec{\agot}^{1},\,
        \vec{f}\,^{-})$ \Comment{Back base change}
        \State
        $\vec{f}\,^{-} \gets f_\ast^{-} \cdot (\lgot_{0} - \lgot_\ast) +
        \sum_{l=0}^{d'} {f}^{1}_l + \vec{f}\,^{-}$ \Comment{Shift vector
          entries}
        \State \Return
        $( \vec{f}\,^{+} + \vec{f}\,^{-} ) \oplus ( \vec{f}\,^{+} -
        \vec{f}\,^{-} )$
        \EndFunction
      \end{algorithmic}
    \end{small}
  \end{linebox}
  \caption{Butterfly Reduction}
  \label{alg:butterfly_reduction}
\end{figure}

\section{Multiplication in the  residue ring of a fiber}\label{sec:resid}

Under the hypotheses and  notation
of Section
\ref{sec:redu} we explain how to multiply in
the residue ring $\LL$ at $B$.        We assume that there exists a $\KK$-rational point $R$
       on $E$ such that $dR\not =O$.
       The residue ring $\LL$ at $B$
       is an algebra of dimension $d$ over
$\KK$. We let $(u_l)_{l\in \ZZ/d\ZZ}$ be the functions on $E$
       defined in Section \ref{sec:ebfu}.
       For every $l$ in $\ZZ/d\ZZ$
       we let \begin{equation}\label{eq:thetal}
         \theta_l = u_l \bmod{B}\end{equation} be the image of  $u_l$ in
       $\LL$. The $\theta_l$ form a basis
       $\Theta$ of $\LL$ over $\KK$. We want to multiply
       two elements
       \[f = \sum_{l\in \ZZ/d\ZZ} f_l\theta_l \text{\,\,\,
       and \,\,\,} g = \sum_{l\in \ZZ/d\ZZ} g_l\theta_l\]  in $\LL$,
       given
       by their coordinates in the  basis $\Theta$.
       We follow~\cite[Section 4.3.4.]{CLEP}.
       We first define  two functions
       \[\cF = \sum_{l\in \ZZ/d\ZZ} f_lu_l \text{\,\,\,
         and \,\,\,}\cG = \sum_{l\in \ZZ/d\ZZ} g_lu_l\]  in $\cL(\gt)$
       such that $f = \cF\bmod{B}$ and $g = \cG\bmod{B}$. As
       a consequence of Equation (\ref{eq:TuAB'}) we
       can decompose
       the product $\cF\, \cG$ as  a sum $\cC+\cD$ where
       \[\cC = \sum_{l\in \ZZ/d\ZZ} (f_l-f_{l-1})(g_l-g_{l-1})\,x_l
       \text{\,\,\, and \,\,\,} \cD=\cF\, \cG - \cC \in \cL(\gt).\]
       Using the method presented in Section \ref{sec:redu} we find
       scalars $C_l$ such that
       \begin{equation*}\label{eq:eqC}
         \cC\equiv \sum_{l\in \ZZ/d\ZZ}C_lu_l  \pmod{B}.\end{equation*}
        It remains to compute the coordinates
       of $\cD$ in the basis $(u_l)_{l\in \ZZ/d\ZZ}$. To this
       end we  evaluate $\cF$ and $\cG$ at the points $R+lt$ for
       $l\in [0,d-1]$ as explained in Section \ref{sec:but}. We also
       evaluate        $\cC$ at the  points $R+lt$. This is achieved in two
       steps. We first reduce $\cC$ modulo $R+\gt$ using the method
       from Section~\ref{sec:redu}. We
       then evaluate
       the resulting function at the $R+lt$ using the method
       in  Section \ref{sec:but}. We  can  thus compute
       \[\cD(R+lt)=\cF (R+lt)\cdot \cG (R+lt)-\cC(R+lt) \text{\,\,\, for \,\,\,}
         l\in [0,d-1].\]
       The   interpolation method from Section
       \ref{sec:int}  finds scalars
       $D_l$ such that
\begin{equation*}\label{eq:eqd}
         \cD =  \sum_{l\in
\ZZ/d\ZZ}D_lu_l.\end{equation*}
        Finally the product of $f$ and $g$
       is
       \[f g =  \sum_{l\in \ZZ/d\ZZ} (C_l+D_l)\theta_l .\]
       Every step in this calculation is achieved in time
       $\cQ d\log (d)$.

       \medskip

       \noindent This finishes the proof of Proposition
        \ref{prop:multlaw} below.
 One can summarize this proof
by saying that the equivariant version of Chudnovsky's multiplication
algorithm
summarized in ~\cite[Lemma  6]{CLEP}
involves convolutions by  constant vectors corresponding
to evaluation, interpolation and reduction. And the main consequence
of the recursive approach presented in Sections \ref{sec:but},
\ref{sec:int}, and \ref{sec:redu} is that these convolution products
are achieved in time $\cQ d \log (d)$ when $d$
is a power of two.

\medskip

The multiplication is depicted in streamlined pseudo-code in
Figure~\ref{alg:butterfly_multiply}.

\begin{figure}[htbp]
  \begin{linebox}[width=0.95\linewidth]{Algorithm}
    \begin{small}
      \begin{algorithmic}
        \Function{NormalBasis\_Multiply}{$\vec{f}, \vec{g}$}
        \State \hspace{0.5cm}$\triangleright$ \parbox[t]{0.85\linewidth}{ The functions
        $\textsc{Butterfly\_Reduce}_R$ and $\textsc{Butterfly\_Reduce}_B$
        denote the reduction routine defined in
        Figure~\ref{alg:butterfly_reduction}, where the constants $\hgot$ and
        $\lgot$ depend on the points $R \in E(\KK)$ and  $b \in E(\Ks)$ as
        specified in Section~\ref{sec:ebff}.}
        \State\State
        \State $\vec{\alpha} \gets \textsc{Butterfly\_Evaluate}\,(\vec{f})$ \State
        $\vec{\beta} \gets \textsc{Butterfly\_Evaluate}\,(\vec{g})$
        \State
        $\vec{H} \gets (\vec{f} - \sigma(\vec{f})) \,\diamond\, (\vec{g} -
        \sigma(\vec{g}))$ \State
        $\vec{h} \gets \textsc{Butterfly\_Reduce}_R\,(\vec{H})$
        \Comment{Reduction modulo $R + \gt$}
        \State $\vec{\gamma} \gets \textsc{Butterfly\_Evaluate}\,(\vec{h})$
        \State $\vec{\delta} \gets \vec{\alpha} \ \diamond\ \vec{\beta} - \vec{\gamma}$
        \State $\vec{k} \gets \textsc{Butterfly\_Interpolate}\,(\vec{\delta})$
        \State $\vec{h} \gets \textsc{Butterfly\_Reduce}_B\,(\vec{H})$
        \Comment{Reduction modulo $b + \gt$ }
        \State \Return $\vec{h} + \vec{k}$ \EndFunction
      \end{algorithmic}
    \end{small}
  \end{linebox}
  \caption{Butterfly Multiplication}
  \label{alg:butterfly_multiply}
\end{figure}

\begin{proposition}[Multiplication]\label{prop:multlaw}
    There exists a constant $\cQ$ such that the following is true.
        Let $\KK$ be a field with odd characteristic.
        Let $E$ be an elliptic curve over $\KK$.
        Let $\delta \geqslant 1$ be an integer.
        Let $t$ be a point of order $d=2^\delta$
        in $E(\KK)$. We assume that there exists  a point $R$
        in $E(\KK)$ such that $dR\not = O$.
         Let   $b$ be
         a point on $E$ defined over a  separable extension
         of $\KK$. Assume   $db\not = O$.
Let
\[B=b+\gt=\{b, b+t, b+2t, \ldots, b+(d-1)t\}.\] Assume  $B$
is left invariant
by the absolute Galois  group of $\KK$.
Let $\LL$ be the residue ring at $B$. Let
$\Theta = (\theta_l)_{l\in \ZZ/l\ZZ}$ be the basis of $\LL$
defined in Equation \ref{eq:thetal}.
        There exists a straight-line program that on input
        $2d$ scalars $(f_l)_{0\leqslant l\leqslant d-1}$
        and $(g_l)_{0\leqslant l\leqslant d-1}$ in $\KK$
        computes
        $d$ scalars $(h_l)_{0\leqslant l\leqslant d-1}$,
        such that
        \[\sum_{l\in \ZZ/d\ZZ} h_l \theta_{l}
          =\left(\sum_{l\in \ZZ/d\ZZ} f_l \theta_{l}  \right)\times
          \left(\sum_{l\in \ZZ/d\ZZ} g_l \theta_{l}  \right)\in \LL.\]
        This straight-line program consists of  no
        more than  $\cQ d\log(d)$
        operations in $\KK$.
\end{proposition}

       We denote by $\otimes$ the multiplication law
       on $\KK^d$ defined in Proposition \ref{prop:multlaw}
       above. We call it an elliptic multiplication law.

\begin{remark}
  In the special case where $b = R$, the pseudo-code in
  Figure~\ref{alg:butterfly_multiply} simplifies significantly. It consists of
  two evaluations,
  \begin{displaymath}
    \vec{\alpha} \gets \textsc{Butterfly\_Evaluate}(\vec{f})
    \quad\text{ and }\quad
    \vec{\beta} \gets \textsc{Butterfly\_Evaluate}(\vec{g})\,,
  \end{displaymath}
followed by a single interpolation,
  \begin{math}
    \textsc{Butterfly\_Interpolate}(\vec{\alpha} \diamond \vec{\beta})\,,
  \end{math}
  which yields the result.
\end{remark}

\section{Implementation and experimental results}\label{sec:impl-exper-results}

A public implementation of the algorithms presented in this work, developed in
the computational algebra system Magma~\cite{magma}, is available
online~\cite{CLGit25}.
It includes the classical Fast Fourier Transform (FFT) using the Cooley--Tukey
algorithm for evaluating or interpolating polynomials in $x$ at $d$-th roots
of unity, as well as our routines for elliptic functions at $d$-torsion
points, both of which require only $O(d \log d)$ operations in
$\mathbb{F}_q$. It also covers evaluating or interpolating polynomials in $x$
at $d$-torsion points of an elliptic curve, which requires $O(d \log^2 d)$
operations using the ECFFT algorithms of~\cite{bckl}.

We present in Table~\ref{tab:timings} the timings for our Magma
routines on a standard laptop, working modulo a 64-bit prime. These results
are provided for comparison purposes only, given the likely overhead of the
Magma interpreter in these highly recursive routines.
The relative differences in timings, namely, a constant factor of about $5$
between Cooley-Tukey and the elliptic butterflies, and asymptotically an
$O(\log d)$ factor for ECFFT are nonetheless significant in practice.
The factor of about 5 reflects, in some sense, the
additional technical cost of working with elliptic functions
versus polynomials in one variable. Equations (\ref{eq:TuAB'})
are a bit more involved than the classical $x^a \cdot x^b=x^{a+b}$.

\begin{table}[htbp]
  \centering
  \begin{tabular}{cccccccccc}
    \toprule
    \multirow{2}{*}{$\log_2 (d)$} && \multicolumn{2}{c}{Cooley-Tukey~\cite{WikiCT}}
    && \multicolumn{2}{c}{Elliptic Butterflies}
    && \multicolumn{2}{c}{ECFFT~\cite{bckl}}
    \\
    \cmidrule(lr){3-4} \cmidrule(lr){6-7} \cmidrule(lr){9-10}
                                && Evaluate & Interpolate && Evaluate (Fig.~\ref{alg:butterfly_evaluate}) & Interpolate (Fig.~\ref{alg:butterfly_interpolate}) && Enter & Exit
    \\
    \midrule
    8  && 0.00 s& 0.01 s&& 0.00 s& 0.01 s&& 0.02 s& 0.04 s
    \\
    9  && 0.01 s& 0.01 s&& 0.02 s& 0.01 s&& 0.05 s& 0.11 s
    \\
    10 && 0.02 s& 0.03 s&& 0.03 s& 0.03 s&& 0.11 s& 0.26 s
    \\
    11 && 0.04 s& 0.03 s&& 0.06 s& 0.07 s&& 0.26 s& 0.62 s
    \\
    12 && 0.03 s& 0.03 s&& 0.13 s& 0.15 s&& 0.60 s& 1.42 s
    \\
    13 && 0.05 s& 0.06 s&& 0.28 s& 0.31 s&& 1.40 s& 3.25 s
    \\
    14 && 0.12 s& 0.12 s&& 0.56 s& 0.66 s&& 3.47 s& 7.47 s
    \\
    15 && 0.25 s& 0.25 s&& 1.17 s& 1.38 s&& 7.65 s& 16.9 s
    \\
    16 && 0.54 s& 0.55 s&& 2.45 s& 2.72 s&& 17.3 s& 38.7 s
    \\
    \bottomrule
  \end{tabular}\medskip

  \caption{Fast evaluation and interpolation timings}
  \label{tab:timings}
\end{table}

           \section{Elliptic bases for finite field extensions}\label{sec:ebff}

      We have shown in~\cite[Theorem 2]{CLEP} how
      evaluation and interpolation of elliptic functions
      provide normal bases for finite field extensions having
      quasi-linear time multiplication. We consider
      in this section extensions of degree $d=2^\delta$,
      a power of two. We construct
      normal  bases allowing  multiplication in time $\cQ d\log (d)$
      for these extensions.
      We prove the theorem below.

      \begin{theorem}[Fast normal bases]\label{th:normalbasis}
        There exists an absolute constant $\cQ$ such that
        the following is true.
        Let $\KK$ be a finite
        field with $q$ elements and odd characteristic. Let
        $\delta \geqslant 1$ be an integer and let $d=2^\delta$.
        Assume that $4d^4\leqslant q$. Let $\LL/\KK$ be a field
        extension of degree $d$. There exists a normal $\KK$-basis
        $\Theta$ of $\LL$ and a straight-line program that
        takes as input the coordinates in $\Theta$
        of two elements in $\LL$ and returns the coordinates
        of their product, at the expense of $\cQ d\log (d)$
        operations in $\KK$.
      \end{theorem}

      \begin{myproof}
      Let $\Ks$ be a separable closure of $\KK$.
Let  $\nu$ be   the $2$-valuation of $q-1$.

\medskip

We first consider the  case $\nu \geqslant \delta+1$.
This is a classical and favorable case because we have enough
roots of unity.
     Let $b$ be an element of order $2^{\delta+\nu}$
          in $\Kss$. Let \[t=b^q/b=b^{q-1}.\] This is an element
     of multiplicative order $d$ in $\KK^*$. So $\KK(b)$
     is a degree $d$  extension of $\KK$. We call it $\LL$.
     We set \[a=b^d.\] This is an element of order $2^\nu$
     in $\KK^*$. So $\LL$
     is a Kummer extension of $\KK$ and
     $(1, b, b^2, \ldots, b^{d-1})$ is a
     basis for $\LL$ over $\KK$. We call it $\Pi$.
      Multiplication of two elements of $\LL$ given
     by their coordinates in the basis  $\Pi$ requires
     $\cQ d\log (d)$ operations in $\KK$ using FFT because we have an element of order
     $2d$ in $\KK^*$. Of course $\Pi$
     is not a normal basis. But the element
     \begin{equation*}\theta = \sum_{l=0}^{d-1}b^l\end{equation*}
       is a normal element in $\LL/\KK$. We denote by $\Theta$
       the associated normal basis and notice that passing from
       $\Theta$ to $\Pi$ is a Fourier transform of order $d$.
       This requires no more than
       $\cQ d\log (d)$ operations in $\KK$ because we have $d$-th roots
       of unity in $\KK^*$ and $d$ is a power of two.

       \medskip

       We now consider the case
       $\nu \leqslant \delta$ and use the constructions in
      ~\cite{CLEP}.  The integer $d_q$ as defined in~\cite[Definition 1]{CLEP}
       is either $d^2$ or $2d^2$. Since $4d^4\leqslant q$ we have
       $d_q\leqslant \sqrt q$ and we can apply~\cite[Lemma 9]{CLEP}.
       There exists an elliptic curve $E$ over $\KK$, a point $t$
       of order $d$ in $E(\KK)$, a point $R$
       in $E(\KK)$ such that $dR\not =O$, and
       a point $b$ in $E(\Ks)$
       such that $db\not =O$ and the conjugate of $b$ by the Frobenius
       of $E/\KK$ is $b+t$. We set $B = b+\gt$. This is a degree $d$
       divisor on $E$. It is  irreducible over $\KK$. The residue field
       at $B$ is a degree $d$ extension field of $\KK$ which we denote by
       $\LL$. The $\theta_l$ introduced in Equation \ref{eq:thetal}
       form a normal basis $\Theta$
       of $\LL$. According to Proposition
       \ref{prop:multlaw} there exists a straight-line program
       that multiplies to elements in $\LL$
       given by their coordinates in the basis
       $\Theta$ at the expense of $\cQ d\log(d)$ operations in
       $\KK$.
       \end{myproof}

       \section{Elliptic Goppa  codes}\label{sec:gop}

       In this section we study a family of Goppa codes
       in genus one allowing particularly fast encoding. We let $p$
       be an odd prime and $q$  a power of $p$.
       Let $d=2^\delta$ with $\delta \geqslant 1$. We set $d'=d/2$ and assume
 that
       \begin{equation}\label{eq:gop}q\geqslant \max(\frac{d^4}{4}, (2d+1)^2+1).\end{equation} The length of the Hasse interval is $4\sqrt q$. So
       there are two consecutive multiples of $d^2$ in it.
       At least one of them is not congruent to $1$ modulo
       $p$.
       Call it $m$.

       Let $c = q+1-m$, then $|c| \le 2\sqrt{q}$ and $m \not\equiv 1 \pmod{p}$
       is equivalent to $c \not\equiv 0 \pmod{p}$.
       By a theorem of Waterhouse~\cite{Waterhouse1969}, for any integer $c$
       with $|c| \le 2\sqrt{q}$ and $p \nmid c$, there exists an ordinary
       elliptic curve $E$ over $\mathbb{F}_q$ with trace $c$, hence with
       exactly $m$ rational points (for prime fields this was already proved
       by Deuring~\cite{Deuring1941}).
       The group of an elliptic curve over a finite field is either cyclic or
       a product of two cyclic groups $C_{n_1} \times C_{n_2}$ with
       $n_1 \mid n_2$. If $d^2$ divides $m$, then $n_2$ must be a multiple of
       $d$. Hence there exists a point of order $d$.

       We deduce that there exists an elliptic curve
       $E$ over $\KK$ such that $E(\KK)$ contains a point $t$
       of order
       $d$. We set \[T=d't.\] Since  $\# E(\KK) \geqslant (\sqrt q -1 )^2 > 4d^2$,  there exists
       a  point $Q$ in $E(\KK)$ such that \[2dQ\not = O.\]
       We denote by
       \[\iota : E\rightarrow E\] the involution
       $P\mapsto -P$ and by
       \[x : E\rightarrow E/\giota\] the quotient map.
       We call \[\cL = \cL (\gt )^\iota\] the subspace of $\cL (\gt)$ fixed by $\iota$.    We can see $\cL$ as the linear space on the genus zero curve $E/\giota$ associated to the divisor
       \[\sum_{l\in [1,d'-1]}[x(lt)].\]
       In particular functions in $\cL$    have no pole at $\{O, T\}$. And  $\cL$ has dimension $d'$.
The functions \[\ell_0=1 \text{\,\,\, and \,\,\,} \ell_l=v_l-v_{-l} \text{\,\,\, for \,\,\,} l\in [1,d'-1]\] form a basis $\ell$ of $\cL$.
       We let
 \[\ev : \cL  \rightarrow \KK^{d}\] be
the evaluation map at the $Q+lt$ for
$0\leqslant l\leqslant d-1$.
Every function $f$ in $\cL$
has degree $\leqslant d-2$ as a function on $E$
and degree $\leqslant d'-1$ as a function
on $E/\giota$. If $f$ vanishes at $d'$ points in $Q+\gt$ then it has
$d'$ zeros on $E/\giota$ also
because $\iota(Q+\gt)$  and $Q+\gt$ do not intersect.
So $f$ must be zero.
We deduce that $\ev$
is injective and its image $C \subset \KK^{d}$ is a linear code
of length $d$, dimension $d'$ and minimum distance $d'+1$.
We can see $C$ as a subcode of a genus one Goppa code or as a genus zero Goppa code.

\medskip

{\bf Encoding}
amounts to  evaluating the map $\ev$ at a function $f$ in $\cL$
given  by its coordinates $(m_l)_{0\leqslant l \leqslant d'-1}$ in the basis $\ell$. To this end we first compute
the coordinates $(n_l)_{l\in \ZZ/\ZZ}$ of $f$ in the basis $v$ of $\cL(\gt)$.  We have
\[n_0=m_0, \, \, n_{d'}=0, \text{\, and \, }n_l=m_l, \text{\, and \, }
  n_{-l}=-m_l \text{\, for \, } 1\leqslant l \leqslant d'-1.\]
We deduce the coordinates of $f$ in the basis $u$ using Equations (\ref{eq:change}) and
(\ref{eq:change2}).
We finally evaluate $f$ at $Q+\gt$ as explained in Section \ref{sec:but}.
So encoding takes time $\cQ d\log (d)$.

\medskip

{\bf Checking} a received word is achieved in time $\cQ d\log (d)$ also
by first interpolating the received values as
explained in Section \ref{sec:int}. We obtain the coordinates $(n_l)_{l\in \ZZ}$ in the basis $v$ of a function
$f\in \cL(\gt)$ taking the received values.  We check that $f\in \cL$ that is
\[n_{d'}=0  \text{\, and \, } n_l+n_{-l}=0 \text{\, for \, } 1\leqslant l \leqslant d'-1.\]
The initial message is given by the coordinates $(m_l)_{0\leqslant l \leqslant d'-1}$
of $f$ in the basis $\ell$ of $\cL$. Namely \[m_0=n_0 \text{\, and \, }
m_l=n_l \text{\, for \, } 1\leqslant l \leqslant d'-1.\]

If the  check is failed we deduce that there are errors. We then remind that $C$ is a genus $0$ code
under its clothes of genus $1$ code. In this context, {\bf decoding} up to half the minimum distance
is achieved thanks to an half-gcd computation as explained in~\cite{SKHN}. The time complexity is $\cQ \log (d)$ times the complexity of multiplying
two polynomials of degree $d$ and coefficients
in $\KK$. See~\cite[Chapter 11]{GG}.
The following theorem summarizes the content of this section

\begin{theorem}[Fast MDS codes]
  There exists a constant $\cQ$ such that the following is true.
  Let $d\geqslant 2$ be a power of two  and let $q$ be an odd prime power
  such that  inequality (\ref{eq:gop}) holds true.
  Let $\KK$ be a field with $q$ elements.
  There exists a $[d,d/2,d/2+1]$ linear code $C$ over $\KK$,
  a straight-line program that encodes $C$ at the expense
  of
  $\cQ d \log (d)$ operations in $\KK$,
  a straight-line program that checks $C$ at the expense
  of
  $\cQ d \log (d)$ operations in $\KK$, and a computation tree that corrects
  $C$ up to  $d/4$ errors at the expense
  of  $\cQ d \log^2 (d)\log (\log (d))$
operations and comparisons in $\KK$.
  \end{theorem}

\section{Elliptic LWE cryptography}\label{sec:ellipt-lwe-crypt}

Our third application concerns the construction of secure cryptographic schemes
within quantum computational models, and in particular in the presence of quantum
algorithms such as those of Shor~\cite{shor1997} and Grover~\cite{Grover1996}.
Among the cryptosystems most extensively studied in this context are those based
on computational problems over Euclidean lattices.\smallskip

The Learning With Errors (LWE) problem, introduced by Regev,
constitutes a foundational hardness assumption underlying these constructions.
\begin{definition}[Short-LWE assumptions~\cite{regev2005}, {\cite[Lemma 4.4]{regev2009}}]\label{def:lwe}
  Let $d, q \geq 2$ be integers, let $\chi$ be an error distribution over $\ZZ$
  (typically a discrete Gaussian), and let $\vec{s} \leftarrow \chi^m$ be a secret
  vector.
  Given pairs $(\mathsf{A}, \mathsf{A}\vec{s} + \vec{e})$, where
  $\mathsf{A} \leftarrow (\ZZ/q\ZZ)^{m \times d}$ is uniformly random and
  $\vec{e} \leftarrow \chi^d$,
  \begin{itemize}
    \item the \emph{short-search-LWE} problem is to recover $\vec{s}$,
    \item the \emph{short-decisional-LWE} problem is to distinguish, with
      non-negligible advantage, such pairs from uniformly random pairs
      $(\mathsf{A}, \vec{u})$.
  \end{itemize}
\end{definition}

Regev establishes that solving the decisional Learning With Errors (LWE)
problem on average is at least as hard as solving worst-case approximation
problems on Euclidean lattices using a quantum algorithm. This remarkable
reduction renders the LWE problem particularly attractive for cryptographic
applications.

We wonder
if replacing matrix products with the elliptic multiplication law defined in
Proposition~\ref{prop:multlaw} could lead to cryptosystems with improved
practical characteristics. This question motivates the definition of an
\emph{Elliptic-LWE} assumption, which transposes LWE into the setting of
elliptic multiplications.

\begin{definition}[Short-Elliptic-LWE assumptions]\label{def:elllwe}
  Let $E$ be an elliptic curve defined modulo a prime $q$, with a point
  $t \in E(\mathbb{Z}/q\mathbb{Z})$ of order $d = 2^\delta$, and let $b$ be
  another point in $E(\mathbb{Z}/q\mathbb{Z})$ such that $d b \neq O$.
   Let
  $\otimes$ denote the multiplication law modulo $b + \langle t \rangle$
  defined in Proposition~\ref{prop:multlaw}. Finally, let $\chi$ be an error
  distribution over $\mathbb{Z}$ (typically a discrete Gaussian), and let
  $\vec{s} \leftarrow \chi^d$ be a secret vector.

  Given pairs $(\vec{a},\, \vec{w} = \vec{a} \otimes \vec{s} + \vec{e})$, where
  $\vec{a} \in (\mathbb{Z}/q\mathbb{Z})^d$ is uniformly random and
  $\vec{e} \leftarrow \chi^d$\,:
  \begin{itemize}
  \item the \emph{short-search Elliptic-LWE} problem is to recover $\vec{s}$,
  \item the \emph{short-decisional Elliptic-LWE} problem is to distinguish such
    pairs from uniformly random pairs $(\vec{a}, \vec{v})$ with non-negligible
    advantage.
  \end{itemize}
\end{definition}

Transposing the Regev construction in this setting yields a secure encryption
scheme in the chosen-plaintext attack (CPA) model.

\begin{theorem}
  Under the Short-Elliptic-LWE assumption, the encryption scheme defined in
  Figure~\ref{fig:ell-pke} is CPA-secure.
\end{theorem}

The correctness follows from the equation
\begin{math}
  p = \langle \vec{r}\,, \vec{e}\rangle - \langle \vec{e}_1\,, \vec{s}\rangle
  + e_2  + \mu\, \lfloor q/2 \rceil \,.
\end{math}
Since the error term remains small for appropriate parameter choices,
decryption succeeds with high probability.  CPA-security follows
directly from the decisional LWE assumption.  The vector $\vec{w}$ is
indistinguishable from a random vector, which makes $\vec{c}_1$ and
$\langle \vec{r},\, \vec{w}\rangle + e_2$ two LWE instances, each
indistinguishable from uniform, thereby perfectly hiding the message.

In terms of complexity, applying  $\phi_{\vec{a}}$ requires
$\mathcal{Q}\, d \log(d)$ operations in $\mathbb{Z}/q\mathbb{Z}$, using the
algorithms described in Section~\ref{sec:redu}.  By Tellegen’s principle,
applying the transpose map $\,^{\mathrm{t}}\phi_{\vec{a}}$ has the same
asymptotic complexity~\cite{KKDB1988,BLS2003}.  More generally, all steps can
be carried out within this time bound.  The sizes of the keys and ciphertexts
are $O(d)$ bits.  \medskip

\begin{figure}[htbp]

  \begin{linebox}[width=0.95\linewidth]{Parameters}
    \begin{small}
      Let $E$ be an elliptic curve defined modulo a prime $q$ with a point
      $t \in E(\mathbb{Z}/q\mathbb{Z})$ of order $d = 2^\delta$, and let $b$
      be another rational point in $E(\mathbb{Z}/q\mathbb{Z})$ such that
      $d \, b \neq O$. Let $\otimes$ denote the multiplication law modulo
      $b + \langle t \rangle$ defined by Proposition~\ref{prop:multlaw}. Let
      $\phi_{\vec{a}}$ (and $\,^\mathrm{t}\phi_{\vec{a}}$\,) denote the linear
      endomorphism $\vec{x}\mapsto \vec{a} \otimes \vec{x}$ in the coordinate
      system defined by the basis $(u_l)_{l \in \ZZ/d\ZZ}$ (and its transposed
      for the canonical  scalar product $\langle \cdot,\, \cdot \rangle$  in this basis).
    \end{small}
  \end{linebox}

  \begin{linebox}[width=0.95\linewidth]{Key Generation}
    \begin{small}
      Sample a uniformly random $d$-dimensional vector $\vec{a}$, a secret
      key $\vec{s} \leftarrow \chi^d$ and an error vector
      $\vec{e} \leftarrow \chi^d$. Compute
      $\vec{w} = \phi_{\vec{a}} ( \vec{s}) + \vec{e}$. The public key is
      $ (\vec{a}, \vec{w})$.
    \end{small}
  \end{linebox}

  \begin{linebox}[width=0.95\linewidth]{Encryption}
    \begin{small}
      To encrypt a message ${\mu} \in \{0,1\}$, sample random
      vectors $\vec{r} \leftarrow  \chi^d$, $\vec{e}_1
      \leftarrow \chi^d$ and $e_2 \leftarrow  \chi$. Then compute
      $\vec{c}_1 = \,^t\phi_{\vec{a}} ( \vec{r}) + \vec{e}_1$ and ${c}_2 =
      \langle \vec{r},\, \vec{w}\rangle + e_2 + \lfloor q/2 \rfloor \,
      {\mu}$. The ciphertext is the pair $(\vec{c}_1, {c}_2)$.
    \end{small}
  \end{linebox}

  \begin{linebox}[width=0.95\linewidth]{Decryption}
    \begin{small}
      To decrypt a ciphertext $(\vec{c}_1, {c}_2)$, compute
      ${p} = {c}_2 - \langle \vec{c}_1, \vec{s}\rangle$.  The
      plaintext bit $\mu$ is then recovered by comparing $p$ to
      $\lfloor q/2 \rceil$: if $p$ is closer to $0$ than to
      $\lfloor q/2 \rceil$ modulo $q$, output $0$; otherwise, output $1$.
    \end{small}
  \end{linebox}

  \caption{Elliptic-LWE encryption scheme variant of Regev's construction}
  \label{fig:ell-pke}
\end{figure}

We note that the elliptic-decisional LWE
assumption bears similarities to the decisional Ring-LWE assumption studied by
Lyubashevsky, Peikert, and Re\-gev~\cite{lyubashevsky2013}, in which matrix
products are replaced by polynomial multiplications in the quotient ring
$\mathbb{Z}_q[x] / \langle \Phi(x) \rangle$ with $\Phi(x)$ a cyclotomic
polynomial.
When the degree $d$ is a power of two, this construction operates in the
cyclotomic field defined by $\Phi_{2d}(x) = x^d + 1$. Since polynomial
multiplication modulo such a polynomial can be implemented using
Fast-Fourier-Transform algorithms provided that $2d$ divides $q-1$, this
approach achieves quasi-linear time complexity in $d$ for encryption and
decryption operations, while also significantly reducing key sizes. This
substantial improvement forms the core of cryptosystems currently undergoing
standardization~\cite{kyber2024,dilithium2024}.

Ideally, one would like in our case to encrypt similarly $d$ bits at once,
as is possible with Ring-LWE. However, the elliptic
multiplication law is not  compatible with the norms on the operands, and
it appears difficult to control the error in the computations; as a result,
the decryption algorithm simply would not work.
The reason for this  incompatibility of the
 elliptic
 multiplication law  with the norms on the operands
 is that the coefficients of the multiplication tensor $\otimes$ have no reason to be small
 as is the case in the cyclotomic context.

 On the other hand, the difficulty of conveniently lifting the Elliptic-LWE
 problem to characteristic zero may be seen as an argument
 for the Short-Elliptic-LWE assumption that the Ring-LWE assumption
 lacks. Also it would be desirable to design  more efficient
 schemes relying on the Short-Elliptic-LWE assumption.

 Nevertheless, it is still possible to make the algorithm in
 Figure~\ref{fig:ell-pke} practical using standard techniques from LWE-based
 cryptography (see, for instance, the Frodo construction~\cite{Frodo2016}).
More precisely, let $\beta$ and $\ell$ be two integer parameters.  We can
encrypt a $\beta\ell^2$-bit message, denoted $(\mu_{i,j})$, by applying
the following modifications to the algorithm in Figure~\ref{fig:ell-pke}.
\begin{itemize}
  \item Set the modulus $q$ so that $\beta$ bits of message can be recovered,
  instead of a single one, in the decryption equation;
  \item Generate secret keys consisting of $\ell$ vectors
  $\vec{s}_0$, $\ldots,\, \vec{s}_{\ell-1}$, associated with public keys
  $(\vec{a}, \vec{w}_0$, $\ldots,\, \vec{w}_{\ell-1})$;
  \item Replace the ciphertext vector $\vec{c}_1$ by $\ell$ vectors computed
  using $\ell$ ephemeral secret vectors
  $\vec{r}_0$, $\ldots,$ $\vec{r}_{\ell-1}$;
\item Compute $\ell \times \ell$ elements $c_2$, one for each inner product
  $\langle \vec{r}_i, \vec{w}_j \rangle_{i,j}$ and addition of
  $e_{i,j} + \lfloor q / 2^{\lfloor 1 + \log_2 q \rfloor  - \beta} \rceil \, \mu_{i,j}$. %
\end{itemize}

The overall time complexity is now $O(\ell\, d \log d + \ell^2 d)$ operations
in $\ZZ/q\ZZ$.  The first term arises from the application of the
endomorphisms $\phi_{\vec{a}}$ and $\,^{\mathrm{t}}\phi_{\vec{a}}$ to the
$\ell$ secret vectors $\vec{s}_i$ and $\vec{r}_j$ during key generation and
encryption.  The second term corresponds to the $\ell^2$ inner products in
the encryption or decryption.  The size of the exchanged data is
$\ell d \log q$ bits.

Using the same parameters as those employed in standards, for instance
$\beta=4$ and $\ell=8$ to encrypt 256-bit data~\cite{Frodo2016}, we observe
that this results in timings and data sizes that are not fundamentally larger
than in the single-bit version.

More generally, most lattice-based cryptographic schemes roughly follow the
following asymptotic guideline, expressed in terms of a security parameter
$\lambda$ tending to infinity~: $d = O(\lambda \log \lambda)$, $q = O(d)$, and
$\chi$ a discrete gaussian distribution with standard deviation
$\sigma = O(\sqrt{d})$~\cite[Chapter 4]{Peikert2016}. We can then extract
$\beta = O(\log d)$ bits per inner product. Setting furthermore
$\ell \simeq O(\sqrt{d }\, / \log d)$ so that we can encrypt on the order of
$\lambda$ bits.
The total time complexity becomes $O(d^{2}\,/\log d)$ bit-operations, while
the data sizes become $O(d^{3/2})$ bits. These complexities are significantly
better than the $O(d^{5/2})$ bit-operations of~\cite{Frodo2016} for same data
size.

\printbibliography

@article {BAL,
  AUTHOR =		 {Ballet, St\'{e}phane},
  TITLE =		 {Curves with many points and multiplication complexity in any
                  extension of {${\bf F}_q$}},
  JOURNAL =		 {Finite Fields Appl.},
  FJOURNAL =	 {Finite Fields and their Applications},
  VOLUME =		 5,
  YEAR =		 1999,
  NUMBER =		 4,
  PAGES =		 {364--377},
  ISSN =		 {1071-5797},
  MRCLASS =		 {11G20 (11R58)},
  MRNUMBER =	 1711833,
  DOI =			 {10.1006/ffta.1999.0255},
}

@book {BCS,
  AUTHOR =		 {B\"urgisser, Peter and Clausen, Michael and Shokrollahi, M.
                  Amin},
  TITLE =		 {Algebraic complexity theory},
  SERIES =		 {Grundlehren der mathematischen Wissenschaften},
  VOLUME =		 315,
  NOTE =		 {With the collaboration of Thomas Lickteig},
  PUBLISHER =	 {Springer-Verlag, Berlin},
  YEAR =		 1997,
  PAGES =		 {xxiv+618},
  ISBN =		 {3-540-60582-7},
  MRCLASS =		 {68-02 (12Y05 65Y20 68Q05 68Q15 68Q25 68Q40)},
  MRNUMBER =	 1440179,
}

@inproceedings{BLS2003,
  author =		 {Bostan, A. and Lecerf, G. and Schost, \'{E}.},
  title =		 {Tellegen's principle into practice},
  year =		 2003,
  isbn =		 1581136412,
  publisher =	 {Association for Computing Machinery},
  address =		 {New York, NY, USA},
  booktitle =	 {Proceedings of the 2003 International Symposium on Symbolic
                  and Algebraic Computation},
  pages =		 {37–44},
  numpages =	 8,
  location =	 {Philadelphia, PA, USA},
  series =		 {ISSAC '03}
}

@article {BR,
  AUTHOR =		 {Ballet, S. and Rolland, R.},
  TITLE =		 {Multiplication algorithm in a finite field and tensor rank
                  of the multiplication},
  JOURNAL =		 {J. Algebra},
  FJOURNAL =	 {Journal of Algebra},
  VOLUME =		 272,
  YEAR =		 2004,
  NUMBER =		 1,
  PAGES =		 {173--185},
  ISSN =		 {0021-8693},
  MRCLASS =		 {11Y16 (11T30)},
  MRNUMBER =	 2029030,
  MRREVIEWER =	 {Igor E. Shparlinski},
  DOI =			 {10.1016/j.jalgebra.2003.09.031},
}

@article{BS,
  author =		 {Bostan, Alin and Schost, {\'E}ric},
  title =		 {Polynomial evaluation and interpolation on special sets of
                  points},
  fjournal =	 {Journal of Complexity},
  journal =		 {J. Complexity},
  issn =		 {0885-064X},
  volume =		 21,
  number =		 4,
  pages =		 {420--446},
  year =		 2005,
  language =	 {English},
  doi =			 {10.1016/j.jco.2004.09.009},
  zbMATH =		 2201759,
  Zbl =			 {1101.68039}
}

@incollection {CHA,
  AUTHOR =		 {Chaumine, Jean},
  TITLE =		 {Multiplication in small finite fields using elliptic curves},
  BOOKTITLE =	 {Algebraic geometry and its applications},
  SERIES =		 {Ser. Number Theory Appl.},
  VOLUME =		 5,
  PAGES =		 {343--350},
  PUBLISHER =	 {World Sci. Publ., Hackensack, NJ},
  YEAR =		 2008,
  MRCLASS =		 {12E20 (14G15)},
  MRNUMBER =	 2484063,
  DOI =			 {10.1142/9789812793430\_0018},
}

@article{CLEP,
  author =		 {Jean-Marc Couveignes and Reynald Lercier},
  title =		 {Elliptic periods for finite fields},
  journal =		 {Finite Fields Their Appl.},
  volume =		 15,
  number =		 1,
  pages =		 {1--22},
  year =		 2009,
}

@software{CLGit25,
  author =		 {Couveignes, Jean-Marc and Lercier, Reynald},
  month =		 oct,
  title =		 {{Elliptic Butterflies}},
  url =			 {https://github.com/rlercier/Elliptic-Butterflies},
  version =		 {1.0.0},
  year =		 2025,
  note =		 {Computer Software},
}

@article {Deuring1941,
    AUTHOR = {Deuring, Max},
     TITLE = {Die {T}ypen der {M}ultiplikatorenringe elliptischer
              {F}unktionenk\"orper},
   JOURNAL = {Abh. Math. Sem. Hansischen Univ.},
  FJOURNAL = {Abhandlungen aus dem Mathematischen Seminar der Hansischen
              Universit\"at},
    VOLUME = {14},
      YEAR = {1941},
     PAGES = {197--272},
      ISSN = {0025-5858},
}

@inproceedings{Frodo2016,
  author =		 {Bos, Joppe and Costello, Craig and Ducas, Leo and Mironov,
                  Ilya and Naehrig, Michael and Nikolaenko, Valeria and
                  Raghunathan, Ananth and Stebila, Douglas},
  title =		 {Frodo: Take off the Ring! Practical, Quantum-Secure Key
                  Exchange from LWE},
  year =		 2016,
  isbn =		 9781450341394,
  publisher =	 {Association for Computing Machinery},
  address =		 {New York, NY, USA},
  booktitle =	 {Proceedings of the 2016 ACM SIGSAC Conference on Computer
                  and Communications Security},
  pages =		 {1006–1018},
  numpages =	 13,
  location =	 {Vienna, Austria},
  series =		 {CCS '16}
}

@book{GG,
  author =		 {Joachim von zur Gathen and J{\"{u}}rgen Gerhard},
  title =		 {Modern Computer Algebra {(3.} ed.)},
  publisher =	 {Cambridge University Press},
  year =		 2013,
  isbn =		 {978-1-107-03903-2},
  bibsource =	 {dblp computer science bibliography, https://dblp.org}
}

@inproceedings{Grover1996,
  author =		 {Grover, Lov K.},
  title =		 {A fast quantum mechanical algorithm for database search},
  year =		 1996,
  isbn =		 0897917855,
  publisher =	 {Association for Computing Machinery},
  address =		 {New York, NY, USA},
  doi =			 {10.1145/237814.237866},
  booktitle =	 {Proceedings of the Twenty-Eighth Annual ACM Symposium on
                  Theory of Computing},
  pages =		 {212–219},
  numpages =	 8,
  location =	 {Philadelphia, Pennsylvania, USA},
  series =		 {STOC '96}
}

@article{HJB,
  author =		 {Heideman, Michael T. and Johnson, Don H. and Burrus,
                  C. Sidney},
  title =		 {Gauss and the history of the fast {Fourier} transform},
  fjournal =	 {Archive for History of Exact Sciences},
  journal =		 {Arch. Hist. Exact Sci.},
  issn =		 {0003-9519},
  volume =		 34,
  pages =		 {265--277},
  year =		 1985,
  language =	 {English},
}

@article {HvdH,
    AUTHOR = {Harvey, David and van der Hoeven, Joris},
     TITLE = {Polynomial multiplication over finite fields in time {$O(n\log
              n)$}},
   JOURNAL = {J. ACM},
  FJOURNAL = {Journal of the ACM},
    VOLUME = {69},
      YEAR = {2022},
    NUMBER = {2},
     PAGES = {Art. 12, 40},
      ISSN = {0004-5411,1557-735X}
}

@article {JV,
  AUTHOR =		 {V{\'e}lu, Jacques},
  TITLE =		 {Isog\'enies entre courbes elliptiques},
  JOURNAL =		 {C. R. Acad. Sci. Paris S\'er. A-B},
  VOLUME =		 273,
  YEAR =		 1971,
  PAGES =		 {A238--A241},
}

@article{JV2,
  author =		 {V{\'e}lu, Jacques},
  title =		 {Courbes elliptiques munies d'un sous-groupe
                  {{\(\mathbb{Z}/n\mathbb{Z}\times\mu_n\)}}},
  fjournal =	 {Bulletin de la Soci{\'e}t{\'e} Math{\'e}matique de
                  France. Suppl{\'e}ment. M{\'e}moires},
  journal =		 {Bull. Soc. Math. Fr., Suppl., M{\'e}m.},
  issn =		 {0583-8665},
  volume =		 57,
  pages =		 152,
  year =		 1978,
  language =	 {French},
  url =			 {https://eudml.org/doc/94779},
  zbMATH =		 3674252,
  Zbl =			 {0433.14029}
}

@article{KKDB1988,
  author =		 {Kaminski, Michael and Kirkpatrick, David G. and Bshouty,
                  Nader H.},
  title =		 {Addition requirements for matrix and transposed matrix
                  products},
  year =		 1988,
  issue_date =	 {September 1988},
  publisher =	 {Academic Press, Inc.},
  address =		 {USA},
  volume =		 9,
  number =		 3,
  issn =		 {0196-6774},
  journal =		 {J. Algorithms},
  month =		 sep,
  pages =		 {354–364},
  numpages =	 11
}

@article{Peikert2016,
  year =		 2016,
  volume =		 10,
  journal =		 {Foundations and Trends® in Theoretical Computer Science},
  title =		 {A Decade of Lattice Cryptography},
  doi =			 {10.1561/0400000074},
  issn =		 {1551-305X},
  number =		 4,
  pages =		 {283-424},
  author =		 {Chris Peikert}
}

@article {RAN,
  AUTHOR =		 {Randriambololona, Hugues},
  TITLE =		 {Bilinear complexity of algebras and the
                  {C}hudnovsky-{C}hud\-novsky interpolation method},
  JOURNAL =		 {J. Complexity},
  FJOURNAL =	 {Journal of Complexity},
  VOLUME =		 28,
  YEAR =		 2012,
  NUMBER =		 4,
  PAGES =		 {489--517},
  ISSN =		 {0885-064X},
  MRCLASS =		 {11G20 (11Y16 14G15 15A72)},
  MRNUMBER =	 2925903,
  MRREVIEWER =	 {Yves Aubry},
  DOI =			 {10.1016/j.jco.2012.02.005},
}

@article {SH,
  AUTHOR =		 {Shokrollahi, Mohammad Amin},
  TITLE =		 {Optimal algorithms for multiplication in certain finite
                  fields using elliptic curves},
  JOURNAL =		 {SIAM J. Comput.},
  FJOURNAL =	 {SIAM Journal on Computing},
  VOLUME =		 21,
  YEAR =		 1992,
  NUMBER =		 6,
  PAGES =		 {1193--1198},
  ISSN =		 {0097-5397},
  MRCLASS =		 {11Y16 (03D15 11G20 68Q25)},
  MRNUMBER =	 1192302,
  MRREVIEWER =	 {Joe P. Buhler},
  DOI =			 {10.1137/0221071},
}

@article{SKHN,
  author =		 {Yasuo Sugiyama and Masao Kasahara and Shigeichi Hirasawa and
                  Toshihiko Namekawa},
  title =		 {A Method for Solving Key Equation for Decoding {G}oppa
                  Codes},
  journal =		 {Inf. Control.},
  volume =		 27,
  number =		 1,
  pages =		 {87--99},
  year =		 1975,
  url =			 {https://doi.org/10.1016/S0019-9958(75)90090-X},
  doi =			 {10.1016/S0019-9958(75)90090-X},
  bibsource =	 {dblp computer science bibliography, https://dblp.org}
}

@incollection {STV,
  AUTHOR =		 {Shparlinski, Igor E. and Tsfasman, Michael A. and Vladut,
                  Serge G.},
  TITLE =		 {Curves with many points and multiplication in finite fields},
  BOOKTITLE =	 {Coding theory and algebraic geometry ({L}uminy, 1991)},
  SERIES =		 {Lecture Notes in Math.},
  VOLUME =		 1518,
  PAGES =		 {145--169},
  PUBLISHER =	 {Springer, Berlin},
  YEAR =		 1992,
  MRCLASS =		 {11G20 (11T71 11Y16)},
  MRNUMBER =	 1186422,
  MRREVIEWER =	 {Jos\'{e} Felipe Voloch},
  DOI =			 {10.1007/BFb0087999},
}

@article {Waterhouse1969,
    AUTHOR = {Waterhouse, William C.},
     TITLE = {Abelian varieties over finite fields},
   JOURNAL = {Ann. Sci. \'Ecole Norm. Sup. (4)},
  FJOURNAL = {Annales Scientifiques de l'\'Ecole Normale Sup\'erieure.
              Quatri\`eme S\'erie},
    VOLUME = {2},
      YEAR = {1969},
     PAGES = {521--560},
      ISSN = {0012-9593}
}

@misc{WikiCT,
  author =		 {Wikipedia},
  title =		 {Cooley--Tukey FFT algorithm},
  url =			 {https://en.wikipedia.org/wiki/Cooley-Tukey_FFT_algorithm},
  note =		 {Online; accessed 2026}
}

@incollection{bckl,
  author =		 {Ben-Sasson, Eli and Carmon, Dan and Kopparty, Swastik and
                  Levit, David},
  title =		 {Elliptic curve fast {Fourier} transform ({ECFFT}). {I}:
                  {Low}-degree extension in time {{\({O}(n\,{{\log}} n)\)}}
                  over all finite fields},
  booktitle =	 {Proceedings of the 34th annual ACM-SIAM symposium on
                  discrete algorithms, SODA 2023},
  isbn =		 {978-1-61197-755-4},
  pages =		 {700--737},
  year =		 2023,
  publisher =	 {New York, NY: Association for Computing Machinery (ACM)},
  language =	 {English},
  doi =			 {10.1137/1.9781611977554.ch30},
  zbMATH =		 7847990
}

@online{wiki,
 author = {Sharp, Charles J. },
 title   = "Pacific double-saddle butterflyfish ({C}haetodon ulietensis) and other {C}hae\-todon, {M}oorea --- Wikimedia Commons",
 year    = "2024",
 urlseen = "03-09-25",
 url     = "https://commons.wikimedia.org/wiki/File:Pacific_double-saddle_butterflyfish_(Chaetodon_ulietensis)_and_other_Chaetodon_Moorea.jpg",
}

@misc{cg,
  title =		 {Explicit {R}iemann-{R}och spaces in the {H}ilbert class
                  field},
  author =		 {Jean-Marc Couveignes and Jean Gasnier},
  year =		 2023,
  eprint =		 {2309.06754},
  archivePrefix ={arXiv},
  primaryClass = {math.NT},
  url =			 {https://arxiv.org/abs/2309.06754},
}

@misc{chchch,
  author =		 {Chudnovsky, D. V. and Chudnovsky, G. V.},
  title =		 {Computational problems in arithmetic of linear differential
                  equations. {Some} diophantine applications},
  year =		 1989,
  language =	 {English},
  howpublished = {Number theory, {Semin}., {New} {York}/{NY} 1985-88,
                  {Lect}. {Notes} {Math}. 1383, 12-49 (1989).},
  zbMATH =		 4127318,
  Zbl =			 {0688.10033}
}

@article{chu,
  author =		 {Chudnovsky, D. V. and Chudnovsky, G. V.},
  title =		 {Algebraic complexities and algebraic curves over finite
                  fields},
  fjournal =	 {Journal of Complexity},
  journal =		 {J. Complexity},
  issn =		 {0885-064X},
  volume =		 4,
  number =		 4,
  pages =		 {285--316},
  year =		 1988,
  language =	 {English},
  doi =			 {10.1016/0885-064X(88)90012-X},
  zbMATH =		 4094802,
  Zbl =			 {0668.68040}
}

@techreport{dilithium2024,
  title =		 {CRYSTALS-Dilithium Algorithm Specifications and Supporting
                  Documentation},
  author =		 {Bai, Shi and others},
  institution =	 {National Institute of Standards and Technology},
  year =		 2024,
  note =		 {NIST PQC Standard}
}

@article {garsti,
  AUTHOR =		 {García, Arnaldo and Stichtenoth, Henning},
  TITLE =		 {A tower of {A}rtin-{S}chreier extensions of function fields
                  attaining the {D}rinfeld-{V}ladut bound},
  JOURNAL =		 {Invent. Math.},
  FJOURNAL =	 {Inventiones Mathematicae},
  VOLUME =		 121,
  YEAR =		 1995,
  NUMBER =		 1,
  PAGES =		 {211--222},
  ISSN =		 {0020-9910,1432-1297},
  MRCLASS =		 {11G20 (11R58 14G15 94B27)},
  MRNUMBER =	 1345289,
  MRREVIEWER =	 {Jos\'e\ Felipe\ Voloch}
}

@article {gop1,
  AUTHOR =		 {Goppa, V. D.},
  TITLE =		 {Codes on algebraic curves},
  JOURNAL =		 {Dokl. Akad. Nauk SSSR},
  FJOURNAL =	 {Doklady Akademii Nauk SSSR},
  VOLUME =		 259,
  YEAR =		 1981,
  NUMBER =		 6,
  PAGES =		 {1289--1290},
  ISSN =		 {0002-3264},
  MRCLASS =		 {94B05 (14H45)},
  MRNUMBER =	 628795,
  MRREVIEWER =	 {E.\ J. F. Primrose},
}

@article {gop2,
  AUTHOR =		 {Goppa, V. D.},
  TITLE =		 {Algebraic-geometric codes},
  JOURNAL =		 {Izv. Akad. Nauk SSSR Ser. Mat.},
  FJOURNAL =	 {Izvestiya Akademii Nauk SSSR. Seriya Matematicheskaya},
  VOLUME =		 46,
  YEAR =		 1982,
  NUMBER =		 4,
  PAGES =		 {762--781, 896},
  ISSN =		 {0373-2436},
  MRCLASS =		 {94B05 (10F45 14C40 94B35)},
  MRNUMBER =	 670165,
  MRREVIEWER =	 {A.\ Peth\H{o}},
}

@article {ihara,
  AUTHOR =		 {Ihara, Yasutaka},
  TITLE =		 {Some remarks on the number of rational points of algebraic
                  curves over finite fields},
  JOURNAL =		 {J. Fac. Sci. Univ. Tokyo Sect. IA Math.},
  FJOURNAL =	 {Journal of the Faculty of Science. University of Tokyo.
                  Section IA. Mathematics},
  VOLUME =		 28,
  YEAR =		 1981,
  NUMBER =		 3,
  PAGES =		 {721--724 (1982)},
  ISSN =		 {0040-8980},
  MRCLASS =		 {14G15 (14G13 14H25 94B05)},
  MRNUMBER =	 656048,
  MRREVIEWER =	 {M.\ A.\ Kenku},
}

@techreport{kyber2024,
  title =		 {CRYSTALS-Kyber Algorithm Specifications and Supporting
                  Documentation},
  author =		 {Avanzi, Roberto and others},
  institution =	 {National Institute of Standards and Technology},
  year =		 2024,
  note =		 {NIST PQC Standard}
}

@inproceedings{laso,
  author =		 {Alexandre Soro and J{\'{e}}r{\^{o}}me Lacan},
  title =		 {{FNT}-Based {R}eed-{S}olomon Erasure Codes},
  booktitle =	 {7th {IEEE} Consumer Communications and Networking
                  Conference, {CCNC} 2010, Las Vegas, NV, USA, January 9-12,
                  2010},
  pages =		 {1--5},
  publisher =	 {{IEEE}},
  year =		 2010,
}

@article{lyubashevsky2013,
  title =		 {On Ideal Lattices and Learning with Errors over Rings},
  author =		 {Lyubashevsky, Vadim and Peikert, Chris and Regev, Oded},
  journal =		 {Journal of the ACM},
  volume =		 60,
  number =		 6,
  pages =		 {43:1--43:35},
  year =		 2013,
  publisher =	 {ACM}
}

@article{magma,
  author =		 {Bosma, Wieb and Cannon, John and Playoust, Catherine},
  title =		 {The {M}agma algebra system. {I}. {T}he user language},
  note =		 {Computational algebra and number theory (London, 1993)},
  journal =		 {J. Symbolic Comput.},
  fjournal =	 {Journal of Symbolic Computation},
  volume =		 24,
  year =		 1997,
  number =		 {3-4},
  pages =		 {235--265},
  issn =		 {0747-7171},
  mrclass =		 {68Q40},
  mrnumber =	 1484478,
  doi =			 {10.1006/jsco.1996.0125},
  url =			 {https://doi.org/10.1006/jsco.1996.0125}
}

@inproceedings{regev2005,
  title =		 {On Lattices, Learning with Errors, Random Linear Codes, and
                  Cryptography},
  author =		 {Regev, Oded},
  booktitle =	 {Proceedings of the 37th Annual ACM Symposium on Theory of
                  Computing},
  pages =		 {84--93},
  year =		 2005
}

@inproceedings{regev2009,
  title =		 {On lattices, learning with errors, random linear codes, and
                  cryptography},
  author =		 {Regev, Oded},
  booktitle =	 {Proceedings of the 41st annual ACM symposium on Theory of
                  computing (STOC)},
  pages =		 {84--93},
  year =		 2009,
  organization = {ACM}
}

@article{shor1997,
  title =		 {Polynomial-Time Algorithms for Prime Factorization and
                  Discrete Logarithms on a Quantum Computer},
  author =		 {Shor, Peter W.},
  journal =		 {SIAM Journal on Computing},
  volume =		 26,
  number =		 5,
  pages =		 {1484--1509},
  year =		 1997,
  publisher =	 {SIAM}
}

@article {tvz,
  AUTHOR =		 {Tsfasman, M. A. and Vl\u{a}du\c{t}, S. G. and Zink, Th.},
  TITLE =		 {Modular curves, {S}himura curves, and {G}oppa codes, better
                  than {V}arshamov-{G}ilbert bound},
  JOURNAL =		 {Math. Nachr.},
  FJOURNAL =	 {Mathematische Nachrichten},
  VOLUME =		 109,
  YEAR =		 1982,
  PAGES =		 {21--28},
  ISSN =		 {0025-584X,1522-2616},
  MRCLASS =		 {11T71 (14G15 14H25 94B05)},
}

\end{document}

%%% Local Variables:
%%% mode: latex
%%% TeX-master: t
%%% End: